\documentclass{article}
\usepackage[legalpaper,margin=1in]{geometry}
\usepackage{graphicx} % Required for inserting images
\usepackage{subcaption}
\usepackage{multirow}
\usepackage[affil-it]{authblk}
\usepackage[hidelinks, breaklinks=true]{hyperref}

\usepackage{amsmath,amsfonts,amssymb}
\usepackage{xcolor}
\usepackage{tikz}
\usetikzlibrary{positioning}
\usetikzlibrary{arrows,calc}
\usepackage{pgfplots}
\newdimen\iwidth
\newdimen\iheight
\newdimen\vsp
\newdimen\hsp
\tikzstyle{every node}=[font=\small]
\usetikzlibrary{spy}
\tikzstyle{label}=[color=white,
fill=darkgray,opacity=.8,
rectangle,rounded corners,
draw,minimum height=0.65cm]
\pgfplotsset{compat=1.18} 

\usepackage{dsfont}
\newcommand{\R}{\ensuremath{\mathds{R}}}

\title{PyHySCO: GPU-Enabled Susceptibility Artifact Distortion Correction in Seconds}
\author[1,*]{Abigail Julian}
\author[1,2]{Lars Ruthotto}
\affil[1]{Department of Computer Science, Emory University, Atlanta, GA, USA}
\affil[2]{Department of Mathematics, Emory University, Atlanta, GA, USA}
\affil[*]{Corresponding Author (abigail.julian@emory.edu)}
\date{}

\begin{document}

\maketitle

\begin{abstract}
    Over the past decade, reversed Gradient Polarity (RGP) methods have become a popular approach for correcting susceptibility artifacts in Echo-Planar Imaging (EPI). Although several post-processing tools for RGP are available, their implementations do not fully leverage recent hardware, algorithmic, and computational advances, leading to correction times of several minutes per image volume. To enable 3D RGP correction in seconds, we introduce PyHySCO, a user-friendly  EPI distortion correction tool implemented in PyTorch that enables multi-threading and efficient use of graphics processing units (GPUs). PyHySCO uses a time-tested physical distortion model and mathematical formulation and is, therefore, reliable without training. An algorithmic improvement in PyHySCO is its novel initialization scheme that uses 1D optimal transport. PyHySCO is published under the GNU public license and can be used from the command line or its Python interface. Our extensive numerical validation using 3T and 7T data from the Human Connectome Project suggests that PyHySCO achieves accuracy comparable to that of leading RGP tools at a fraction of the cost. We also validate the new initialization scheme, compare different optimization algorithms, and test the algorithm on different hardware and arithmetic precision.  
\end{abstract}


\section{Introduction}
\begin{figure}[t]
    \centering
    \input{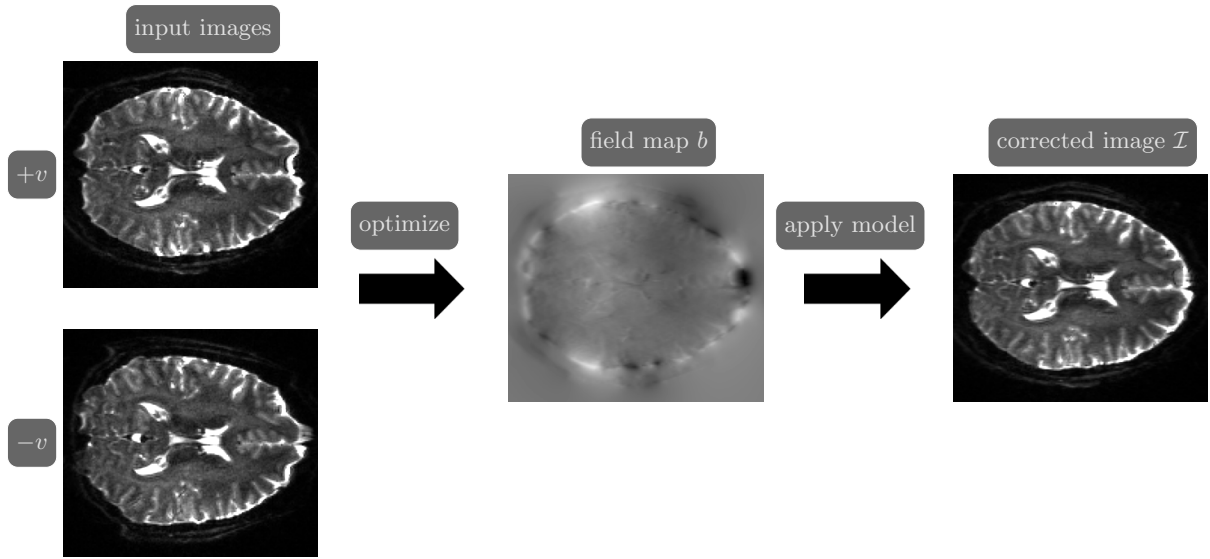}
    \caption{The Reverse Gradient Polarity correction paradigm. Two images are acquired with opposite phase encoding directions, $+v$ and $-v$. These two images are used to estimate the field map $b$, and the distortion correction model \cite{ChangFitzpatrick1992} is applied to obtain a corrected image $\mathcal{I}$.}
    \label{fig:rgp}
\end{figure}

Reversed Gradient Polarity (RGP) methods are commonly used to correct susceptibility artifacts in Echo-Planar Imaging (EPI) \cite{Stehling1991EPI}. RGP methods acquire a pair of images with opposite phase encoding directions, which leads to a minimal increase in scan time due to the speed of EPI. In a post-processing step, RGP approaches use the fact that  the distortion in both images has an equal magnitude but acts in opposite directions to estimate the field map (see Figure \ref{fig:rgp}) \cite{ChangFitzpatrick1992, bowtell1994correction}. The field map is then used to estimate a distortion-free image. 

Compared to other correction approaches such as field map acquisition, point-spread function map acquisition, and anatomical registration, RGP methods generally achieve comparable or superior accuracy while being more robust to  noise and motion; see, e.g., \cite{tax2022s, graham2017quantitative, esteban2014simulation, wu2008comparison}. These advantages make RGP correction a popular choice. For example, the widely-used MRI database from the Human Connectome Project (HCP) \cite{HCPdata} used the RGP correction tool TOPUP \cite{AnderssonEtAl2003} in the preprocessing of released diffusion MRI from EPI scans.

The original RGP distortion correction approaches in \cite{ChangFitzpatrick1992, bowtell1994correction} are one-dimensional, treating each image column separately in the phase encoding direction. This leads to a non-smooth field map estimate and corrections. TOPUP addresses this non-smoothness with a 3D spline-based approach and the introduction of regularization \cite{AnderssonEtAl2003}. TOPUP has limited support for hyperthreading and is often a time-consuming step of MRI processing pipelines \cite{cai2021prequal}. Running TOPUP on a standard CPU took over 60 minutes on average per HCP subject in our experiments. 

Although less widely used than TOPUP, other iterative methods have proposed implementations of RGP correction employing various optimization schemes, discretizations, and regularization terms to speed up the correction. EPIC \cite{HollandEtAl2010} introduces correction using a nonlinear image registration framework. The tool was developed specifically for Anterior-Posterior distortions and can be less effective for Left-Right distortions \cite{gu2019evaluation}. DR-BUDDI \cite{irfanoglu2015dr} and TISAC \cite{duong2020susceptibility} regularize the optimization using a T2-weighted or T1-weighted image, respectively. Including undistorted anatomical information can improve the quality of distortion correction \cite{gu2019evaluation}, but complicates the choice of an effective distance measure and, depending on the protocol, may require additional scan time. HySCO introduces hyper-elastic registration regularization and a novel separable discretization \cite{RuthottoEtAl2012EPI, RuthottoEtAl2013hysco, MacdonaldRuthotto2017}. HySCO can accurately correct real and simulated data varying in phase encoding direction, anatomy, and field of view \cite{gu2019evaluation, Snoussi2021_spinalcordeval, tax2022s}. In our experiments, on average, HySCO runs on the CPU for one to two minutes per HCP subject. While HySCO is a Statistical Parametric Mapping (SPM)~\cite{SPM2007} plugin and has been integrated into a comprehensive SPM-based DTI processing pipeline~\cite{Dvid2024ACIDAC}, its dependency on a MATLAB license may limit its wider use.

Motivated by the long processing times of the above RGP tools, several deep learning approaches for susceptibility artifact correction have been proposed recently; see, e.g., \cite{Hu2020, duong2020unsupervised, Duong2021, zahneisen2020deep, alkilani2023fd}. A recurrent theme is to train a correction operator in an offline stage in a supervised way using training data, which enables fast evaluations in the online step. For example, training S-Net on 150 volumes took over 5 days, while correcting an image pair on a CPU took an average of 2.8 seconds (0.96 seconds on a GPU) \cite{duong2020unsupervised}. 
However, the dramatic reduction of correction time comes at the cost of losing the robustness and generalizability that existing RGP approaches obtain from the physical distortion model. For example, while RGP approaches can handle images from different scanners, anatomies, resolutions, and other acquisition parameters, deep learning models perform poorly when applied outside the training distribution, \cite{chen2022DLReview}. Furthermore, deep learning models are highly sensitive to noise and adversarial attacks in other contexts~\cite{AntunEtAl2020}.

The PyHySCO toolbox aims to achieve the accuracy, robustness, and generalizability of state-of-the-art RGP approaches at computational costs similar to evaluating a pre-trained deep learning model.  PyHySCO offers EPI distortion correction through a GPU-enabled and command-line-accessible Python tool powered by PyTorch~\cite{pytorch}. The mathematical formulation is based on HySCO augmented by the separable discretization  \cite{MacdonaldRuthotto2017} that increases parallelism and an optimal transport-based initialization step that alleviates the need for multilevel optimization.   We demonstrate the use of PyHySCO using its Python interface and command-line tool, which is compatible with existing MRI postprocessing pipelines.

The remainder of our paper is organized as follows.
In Section \ref{sec:2methods}, we review the mathematical model and its discretization under the hood of  PyHySCO and introduce our novel parallelized initialization using optimal transport, fast solvers exploiting the separable structure, and GPU-enabled PyTorch implementation. In Section \ref{sec:3results}, we extensively validate PyHySCO on real and simulated EPI data. We show the speed and accuracy of the proposed parallelized optimal transport initialization scheme and the speed and accuracy of the complete correction pipeline across optimizers, GPUs, and levels of numerical precision. In Section \ref{sec:4discussion}, we discuss the benefits and implications of using PyHySCO for EPI distortion correction. In Section \ref{sec:conclusion}, we provide a conclusion.

\section{Methods}\label{sec:2methods}
This section describes the algorithmic and coding structure of PyHySCO.
Section~\ref{sub:math} introduces our notation and reviews the mathematical formulation of the RGP correction problem. Section~\ref{sub:OTinit} introduces a novel optimal transport-based field map initialization. Section~\ref{sub:optim} describes the optimization algorithms available in PyHySCO. Section~\ref{sec:implementation} explains the structure of the code and some key implementation details. Section~\ref{sub:usage} demonstrates the basic usage of PyHySCO and how to integrate it into existing processing pipelines.

\subsection{Mathematical Formulation}
\label{sub:math}
The field map estimation and distortion correction are based on the physical forward model defined in \cite{ChangFitzpatrick1992}. Let $v \in \R^3$ be the phase encoding direction for the distorted observation $\mathcal{I} : \Omega \rightarrow \R$, and let $\Omega \subset \R^3$ be the image domain of interest. The mass-preserving transformation operator that, given the field map $b : \Omega \rightarrow \R$,  corrects the distortions of an image $\mathcal{I}$ acquired with phase-encoding direction $v$ reads
\begin{equation}\label{eq:transform_model}
    \mathcal{T}[\mathcal{I}, b, v](x)= \mathcal{I}(x+b(x)v) \cdot (1+\partial_{v}b)(x).
\end{equation}
Here, $\partial_{v}b$ is the directional derivative of $b$ in the direction of $v$. The first term of the operator corrects the geometric deformation in the direction of $v$, and the second is an intensity modulation term, which should always be positive.

Similar to \cite{RuthottoEtAl2013hysco}, PyHySCO solves the inverse problem of estimating the field map $b$ based on two observations, $\mathcal{I}_{+v}$ and $\mathcal{I}_{-v}$, acquired with phase-encoding directions $\pm v$, respectively. To this end, we estimate the field map $b$ by minimizing the distance of the corrected images 
\begin{equation}
\label{eq:dist}
    \mathcal{D}(b) = \frac{1}{2}\int_{\Omega} \left(\mathcal{T}[\mathcal{I}_{+v}, b, v](x) - \mathcal{T}[\mathcal{I}_{-v}, b, -v](x)\right)^2 \ \textrm{d}x.
\end{equation}
The distance term is additionally regularized to enforce smoothness and the intensity modulation constraint. The smoothness regularization term, 
\begin{equation*}
    \mathcal{S}(b) = \frac{1}{2}\int_{\Omega}||\nabla b(x)||^2 \textrm{d}x,
\end{equation*}
penalizes large values of the gradient of $b$ to ensure smoothness in all directions.

The intensity modulation constraint of the physical model requires that $-1 < \partial_{v}b(x) < 1$ for almost all $x \in \Omega$. This is enforced by the barrier term
\begin{equation}\label{eq:phi}
    \mathcal{P}(b) = \frac{1}{2}\int_{\Omega} \phi(\partial_{v}b(x)) \textrm{d}x, \ \text{where} \ \phi(z) = \begin{cases}
        \frac{z^4}{1-z^2}, & z \in (-1,1)\\
        \infty, & {\text{else.}}
    \end{cases} 
\end{equation}
Altogether, this gives the optimization problem
\begin{equation}\label{eq:opt}
    \min_{b} \mathcal{J}(b) = \mathcal{D}(b) + \alpha \mathcal{S}(b) + \beta \mathcal{P}(b),
\end{equation}
where the importance of the regularization terms is weighted with non-negative scalars $\alpha$ and $\beta$. Higher values of $\alpha$ promote a smoother field map, while lower values of $\alpha$ promote reduced distance between corrected images at the expense of smoothness in the field map. Any positive value for $\beta$ ensures the intensity modulation constraint is satisfied, but lower values can lead to more ill-conditioned problems. For this paper, we fix $\alpha=300$ and $\beta=1\text{e}-4$.

PyHySCO follows the discretize-then-optimize paradigm commonly used in image registration; see, e.g.,~\cite{Modersitzki2009}.
PyHySCO discretizes the variational problem \eqref{eq:opt} as in \cite{MacdonaldRuthotto2017} to obtain a finite-dimensional optimization problem almost entirely separable in the phase encoding direction. Specifically, coupling is only introduced in the smoothness regularization term when calculating the gradient in the frequency encoding and slice selection directions. 

Our convention is to permute the dimensions of the input image such that the phase encoding direction is aligned with the third unit vector $e_3 = [0,0,1]^T$. The field map is discretized on an $e_3$-staggered grid; that is, we discretize its values in the cell centers along the first two dimensions and on the nodes in the third dimension. The integrals in \eqref{eq:opt} are approximated by a midpoint quadrature rule. The input images are modeled by a one-dimensional piecewise linear interpolation function in the phase encoding direction.The geometric transformation is estimated in the cell centers with an averaging operator, and the intensity modulation is estimated in the cell centers with a finite difference operator. 

The discretized smoothness regularization term is computed for the discretized field map $\mathbf{b}$ via
\begin{equation}\label{eq:smoothness}
    S(\mathbf{b}) = \frac{h_1 \cdot h_2 \cdot h_3}{2} \mathbf{b}^\top H \mathbf{b} = \frac{h_1 \cdot h_2 \cdot h_3}{2}||\mathbf{b}||_H^2,
\end{equation}
where $h_1, h_2, h_3$ are the voxel sizes and $H$ is a standard five-point discretization of the negative Laplacian and thus is a positive semi-definite operator. The discretized intensity modulation constraint term applies $\phi$ as defined in \eqref{eq:phi} element-wise to the result of a finite difference operator applied to the discretized field map. This gives the discretized optimization problem to solve as
\begin{equation}\label{eq:opt_disc}
    \min_{\mathbf{b}} J(\mathbf{b}) = D(\mathbf{b}) + \alpha S(\mathbf{b}) + \beta P(\mathbf{b}).
\end{equation}
This problem is challenging to solve because it is high-dimensional and non-convex, but we can exploit the structure and separability to efficiently solve the problem using parallelization. The implementation of this optimization problem in a parallelizable way, as described in Section \ref{sec:implementation}, includes choices of image interpolation, linear operators for averaging and finite difference, and regularization terms $S$ and $P$.

\subsection{Parallelized Initialization using Optimal Transport}
\label{sub:OTinit}
Due to the non-convexity of the optimization problem \eqref{eq:opt_disc}, an effective initialization strategy for the field map is critical. To this end, PyHySCO introduces a novel initialization scheme based on techniques from optimal transport (OT) \cite{peyreOT}. 
The key idea is to compute the 'halfway' point of the oppositely distorted images in the Wasserstein space (as opposed to Euclidean space, which would simply average the images). 
To render this problem feasible, we treat each image column separately, use the closed-form solutions of 1D OT problems, and then apply a smoothing filter.

We calculate these transformations as optimal transport maps \cite{peyreOT}. More specifically, because the distortions only occur in the phase encoding direction, these transformations are a set of one-dimensional maps calculated in parallel across the distortion dimension. One-dimensional optimal transport has a closed-form solution arising from considering the one-dimensional signal as a positive measure and constructing a cumulative distribution function \cite{peyreOT}. 

We describe the computation of the one-dimensional optimal transport maps in the distortion correction setting. In practice, the computation is parallelized in the distortion dimension to compute the entire initial field map simultaneously. 

Let $i_{+v} \in \R^m$ be the image data from an entry in the phase encoding dimension of $\mathcal{I}_{+v}$, and let $i_{-v} \in \R^m$ be the image data from the corresponding entry in the phase encoding dimension of $\mathcal{I}_{-v}$. Consider $i_{\text{half}}$ the sequence of image intensity values from the corresponding entry of the undistorted image $\mathcal{I}$. We numerically ensure $i_{+v}$ and $i_{-v}$ can be considered positive measures by applying a small shift to the image values, which does not change the relative distance between elements. 

We initialize the field map using the optimal transport maps $T_{+}$ from $i_{+v}$ to $i_{\text{half}}$ and $T_{-}$ from $i_{-v}$ to $i_{\text{half}}$. These maps can be directly computed using the closed-form one-dimensional optimal transport formula, which depends on a cumulative distribution function and its pseudoinverse \cite{peyreOT}.

Define the discretized cumulative distribution function $C_{i}: \{0,\dots,m\} \rightarrow [0,1]$ of a measure $i$ as the cumulative sum
\begin{equation*}
    \forall x \in \{0,\dots,m\} \ \ C_{i}(x) = \sum_{j=0}^x i(j),
\end{equation*}
where $i(j)$ returns the pixel intensity value at index $j$ of $i$. The pseudoinverse $C^{-1}_{i}: [0,1] \rightarrow \{0,\dots,m\}$ is defined as
\begin{equation*}
    \forall r \in [0,1] \ \ C^{-1}_{i}(r) = \min_x \{x \in \{0,\dots,m\} \ | \ C_{i}(x) \geq r\}.
\end{equation*}
In practice, $C^{-1}_{i}$ is computed using a linear spline interpolation.

\begin{figure}[t]
     \centering
    \input {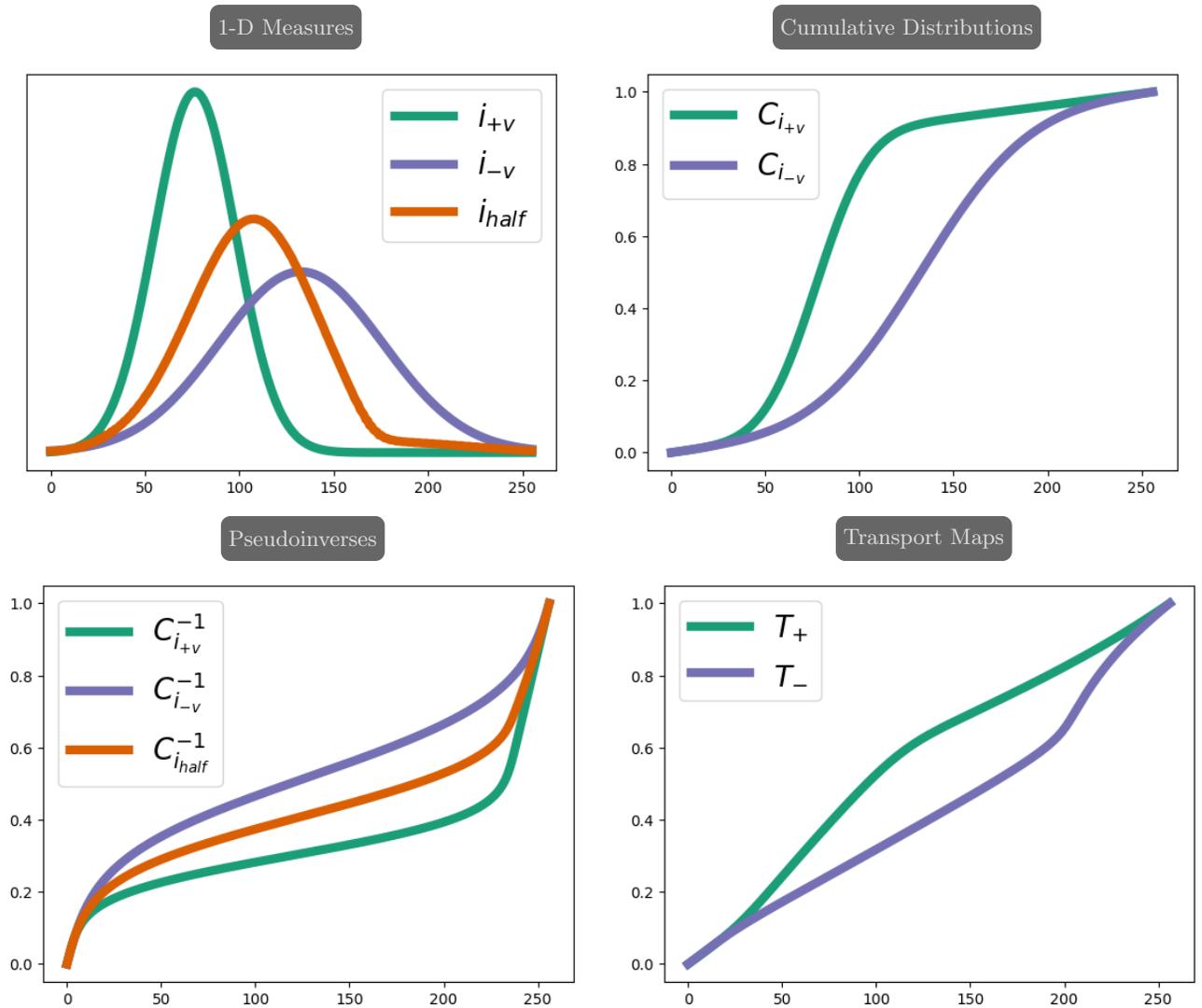}
        \caption{Example of one-dimensional optimal transport maps. The top left shows an example of one-dimensional measures. The green signal, $i_{+v}$, corresponds to an intensity pileup in $\mathcal{I}_{+v}$, while the purple signal $i_{-v}$ corresponds to an intensity dispersion in $\mathcal{I}_{-v}$. The red signal corresponds to the intensity of the true image. The top right shows the cumulative distributions for the measures $i_{+v}$ and $i_{-v}$. Bottom left shows the pseudoinverses for $i_{+v}$ and $i_{-v}$ along with the pseudoinverse $C_{i_{\text{half}}}^{-1}$ used in calculating the transport maps $T_{+}=C^{-1}_{i_{\text{half}}} \circ C_{i_{+v}}$ and $T_{-}=C^{-1}_{i_{\text{half}}} \circ C_{i_{-v}}$, shown bottom right.}
        \label{fig:OT}
\end{figure}

\begin{figure}[t]
    \centering
    \begin{subfigure}{0.9\textwidth}
        \centering
        \input {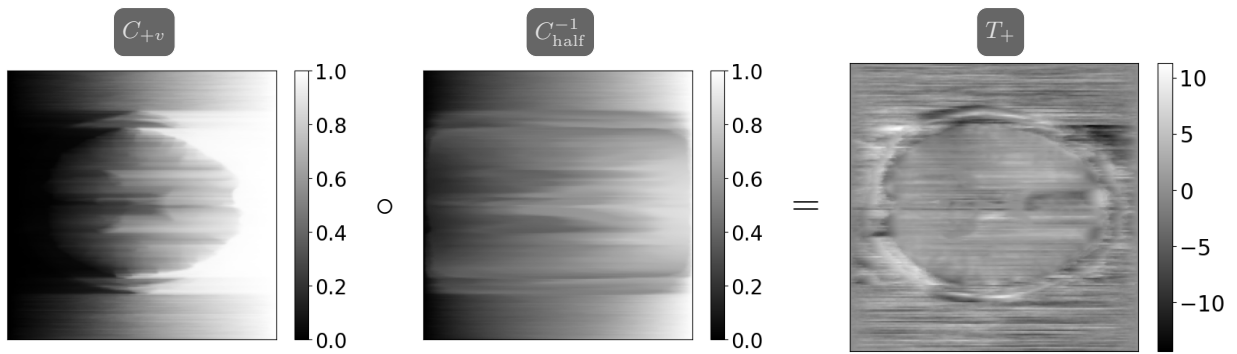}
        \caption{The map $T_{+}$ mapping from $\mathcal{I}_{+v}$ halfway to $\mathcal{I}_{-v}$ is calculated as the composition of the cumulative distribution function $C_{+v}$ from $\mathcal{I}_{+v}$ and the interpolated pseudoinverse $C^{-1}_{\text{half}}$.}
    \end{subfigure}
    \vspace{1em}
    \begin{subfigure}{0.9\textwidth}
        \centering
        \input {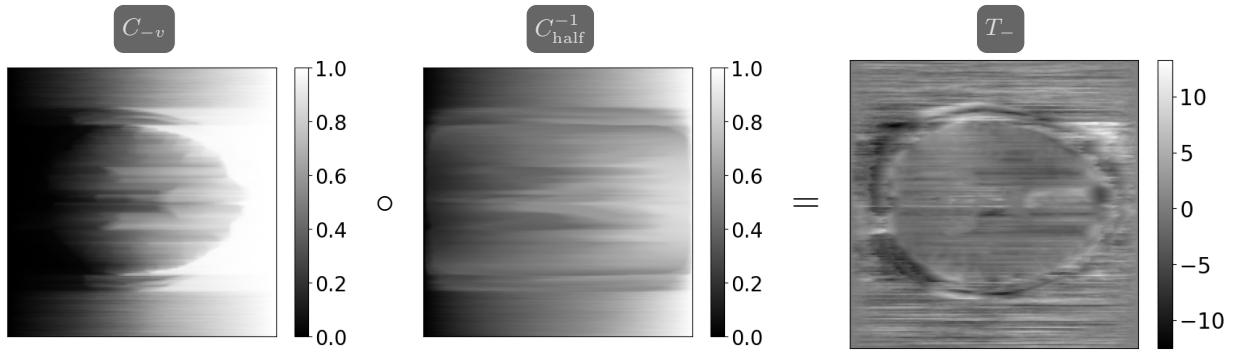}
        \caption{The map $T_{-}$ mapping from $\mathcal{I}_{-v}$ halfway to $\mathcal{I}_{+v}$ is calculated as the composition of the cumulative distribution function $C_{-v}$ from $\mathcal{I}_{-v}$ and the interpolated pseudoinverse $C^{-1}_{\text{half}}$.}
    \end{subfigure}
    \caption{The maps $T_{+}$ and $T_{-}$ are calculated using the closed-form one-dimensional optimal transport solution, parallelized in the distortion dimension \cite{peyreOT}. Note the inverted coloring between $T_{+}$ and $T_{-}$ as the map $T_{-}$ corrects a distortion in the opposite direction as $T_{+}$.}
    \label{fig:ot_maps}
\end{figure}

Returning to the measures arising from the input images, the closed-form solution for one-dimensional optimal transport gives the optimal transport map from $i_{+v}$ to $i_{\text{half}}$ as
\begin{equation*}
    T_{+} = C^{-1}_{i_{\text{half}}} \circ C_{i_{+v}},
\end{equation*}
and the optimal transport map from $i_{-v}$ to $i_{\text{half}}$ as
\begin{equation*}
    T_{-} = C^{-1}_{i_{\text{half}}} \circ C_{i_{-v}},
\end{equation*}
where $C^{-1}_{i_{\text{half}}}$ is calculated as $=(C^{-1}_{i_{+v}}+C^{-1}_{i_{-v}}) / 2$. Figure \ref{fig:OT} visualizes the computation of the one-dimensional transport maps, and the parallelized computation and resulting field maps are visualized in Figure \ref{fig:ot_maps}. We thus compute the initial guess for the field map as the average of the maps $T_{+}$ and $-T_{-}$, computed in parallel. To introduce smoothness in the field map in the frequency encoding and slice selection dimensions, we apply a smoothing filter to the initial field map before optimization.

\subsection{Optimization Algorithms}
\label{sub:optim}
Since the optimal choice of optimization algorithms for approximately solving~\eqref{eq:opt_disc} may depend on various factors, including image sizes, computational hardware, and severity of distortions, PyHySCO offers three options.
Section~\ref{subs:GN} describes a Gauss-Newton scheme with Jacobi-Preconditioned Conjugate Gradient (GN-PCG) method as an inner solver, which is similar to~\cite{RuthottoEtAl2013hysco} and is the default option.
An option that exploits the parallelism of the discretization more effectively is the Alternating Direction Method of Multipliers (ADMM) in Section~\ref{subs:ADMM}, which is based on~\cite{MacdonaldRuthotto2017}.
For comparison, we also provide an interface to an LBFGS optimizer; see Section~\ref{subs:lBFGS}.

\subsubsection{GN-PCG: Gauss-Newton with Jacobi-Preconditioned Conjugate Gradient Solver}
\label{subs:GN}
PyHySCO's default solver is a PyTorch implementation of the GN-PCG scheme used in \cite{RuthottoEtAl2013hysco}.
The key idea is to linearize the residual in the distance term~\eqref{eq:dist} and the nonlinear penalizer~\eqref{eq:phi} about the $k$-th iterate $\mathbf{b}_k$ and approximately solve the resulting quadratic problem with a few iterations of the PCG method. 

More precisely, let $\nabla J$ be the gradient and $H_J$ be a positive definite approximation of the Hessian of the optimization problem \eqref{eq:opt_disc} about $\mathbf{b}_k$. Gauss Newton iteratively updates the current field map estimate via
\begin{equation*}
    \mathbf{b}_{k+1} = \mathbf{b}_k + \gamma_k \textbf{q}_k,
\end{equation*}
where the step size $\gamma_k$ is determined with a line search method such as Armijo \cite[Ch. 3 p. 33-36]{nocedal1999numerical} and the search direction $\textbf{q}_k$ approximately satisfies
\begin{equation}\label{eq:gn_system}
    H_J(\mathbf{b}_k)\textbf{q}_k = - \nabla J(\mathbf{b}_k).
\end{equation}
To obtain $\mathbf{q}_k$, we apply up to 10 iterations of the preconditioned conjugate gradient (PCG) method and stop early if the relative residual is less than $0.1$; see the original work \cite{HestenesStiefel1952} or the textbook \cite{Saad2003} for more details on PCG.
The performance of PCG crucially depends on the clustering of the eigenvalues, which a suitable preconditioner can often improve.
As a computationally inexpensive and often effective option, we implement a Jacobi preconditioner, which approximates the inverse of $H_J$ by the inverse of its diagonal entries. 
Instead of constructing the matrix $H_J$, which is computationally expensive, we provide efficient algorithms to compute matrix-vector products and extract its diagonal.
While the diagonal preconditioner works well in our examples, we note that a more accurate (yet also more expensive) block-diagonal preconditioner has been proposed in  \cite{MacdonaldRuthotto2017}.

\subsubsection{Alternating Direction Method of Multipliers (ADMM)}
\label{subs:ADMM}
We additionally modify the ADMM \cite{boyd2011admm} algorithm in \cite{MacdonaldRuthotto2017} and implement it in PyHySCO.
To take advantage of the separability of the objective function, the idea is to split the optimization problem into two subproblems.
In contrast to \cite{MacdonaldRuthotto2017}, which uses a hard constraint to ensure positivity of the intensity modulation and employs Sequential Quadratic Programming, we implement this as a soft constraint with the barrier term~\eqref{eq:phi}.

As in~\cite{MacdonaldRuthotto2017}, we split the objective in \eqref{eq:opt_disc} into
\begin{align}\label{eq:split}
    F(\mathbf{b}) = D(\mathbf{b}) + \alpha S_3(\mathbf{b}) + \beta P(\mathbf{b}), \quad \text{ and } \quad
    G(\mathbf{z}) = \alpha S_{1}(\mathbf{z}) + \alpha S_{2}(\mathbf{z}),
\end{align}
where $S_3$ is the part of the smoothness regularization term $S$ corresponding to the phase encoding direction, and $S_1$ and $S_2$ are the remaining terms corresponding to the other directions. This gives rise to the following optimization problem, equivalent to \eqref{eq:opt_disc}:
\begin{equation*}
    \min_{\mathbf{b}, \mathbf{z}} F(\mathbf{b}) + G(\mathbf{z}) \ \ \text{s.t.} \ \ \mathbf{b} = \mathbf{z}.
\end{equation*}
With the corresponding augmented Lagrangian
\begin{equation*}
    L(\mathbf{b},\mathbf{z},\mathbf{y}) = F(\mathbf{b}) + G(\mathbf{z}) + \mathbf{y}^T(\mathbf{b}-\mathbf{z}) + \frac{\rho h^3}{2}||\mathbf{b}-\mathbf{z}||^2,
\end{equation*}
where $\mathbf{y}$ is the Langrange multiplier for the equality constraint $\mathbf{b}=\mathbf{z}$ and $\rho$ is a scalar augmentation parameter, and using scaled Lagrange multiplier $\mathbf{u}=\frac{\mathbf{y}}{\rho h^3}$, each iteration has the updates
\begin{align}
    \mathbf{b}_{k+1} &= \arg\min_{\mathbf{b}} F(\mathbf{b}) + \frac{\rho h^3}{2}||\mathbf{b}-\mathbf{z}_k+\mathbf{u}_k||^2 \label{eq:x_update}\\
    \mathbf{z}_{k+1} &= \arg\min_{\mathbf{z}} G(\mathbf{z}) + \frac{\rho h^3}{2}||\mathbf{b}_{k+1}-\mathbf{z}+\mathbf{u}_k||^2 \label{eq:z_update}\\
    \mathbf{u}_{k+1} &= \mathbf{u}_k + \mathbf{b}_{k+1} - \mathbf{z}_{k+1}. \label{eq:u_update}
\end{align}

The $\mathbf{b}$ update computed in~\eqref{eq:x_update} involves a separable optimization problem that can be solved independently for each image column along the phase-encoding direction.
In PyHySCO we use a modified version of the GN-PCG scheme described above.
The only change is the computation of the search direction, which can now be parallelized across the different image columns.
To exploit this structure, we implement a PCG method that solves the system for each image column in parallel. 
In addition to more parallelism, we observe an increase in efficiency since the scheme uses different step sizes and stopping criteria for each image column.

The $\mathbf{z}$ update is computed by solving the quadratic problem \eqref{eq:z_update} directly. This is enabled by the structure of the associated linear system, which is block-diagonal, and each block is given by a 2D negative Laplacian (from the regularizers) shifted by an identity (from the proximal term). 
Assuming periodic boundary conditions on the images, the blocks in the approximation itself have an exploitable structure (called Block Circulant - Circulant Block in~\cite{hansen2006deblurring}) and, therefore,  can be inverted efficiently with the Fast Fourier Transform (FFT). 

The augmentation parameter $\rho$ is updated adaptively as described in \cite{boyd2011admm} to keep the relative primal and dual residuals close.

\subsubsection{LBFGS}
\label{subs:lBFGS}
As a comparison, we provide an implementation of LBFGS \cite{lbfgs}, although optimization with LBFGS does not exploit any of the structure or separability of the optimization problem. LBFGS is a quasi-Newton method that uses an estimate of the objective function's Hessian based on a limited number of previous iterations in solving for the search direction \cite{lbfgs}. In our implementation, we provide an explicitly calculated derivative to an LBFGS solver\footnote{https://github.com/hjmshi/PyTorch-LBFGS}. In computing the objective function, we precompute parts of the derivative which allows for faster optimization than relying on automatic differentiation.

\subsection{Coding Structure of PyHySCO}\label{sec:implementation}
\begin{figure}
    \centering
    \begin{subfigure}{0.9\textwidth}
        \includegraphics[width=\textwidth, trim=175 0 240 0, clip=true]{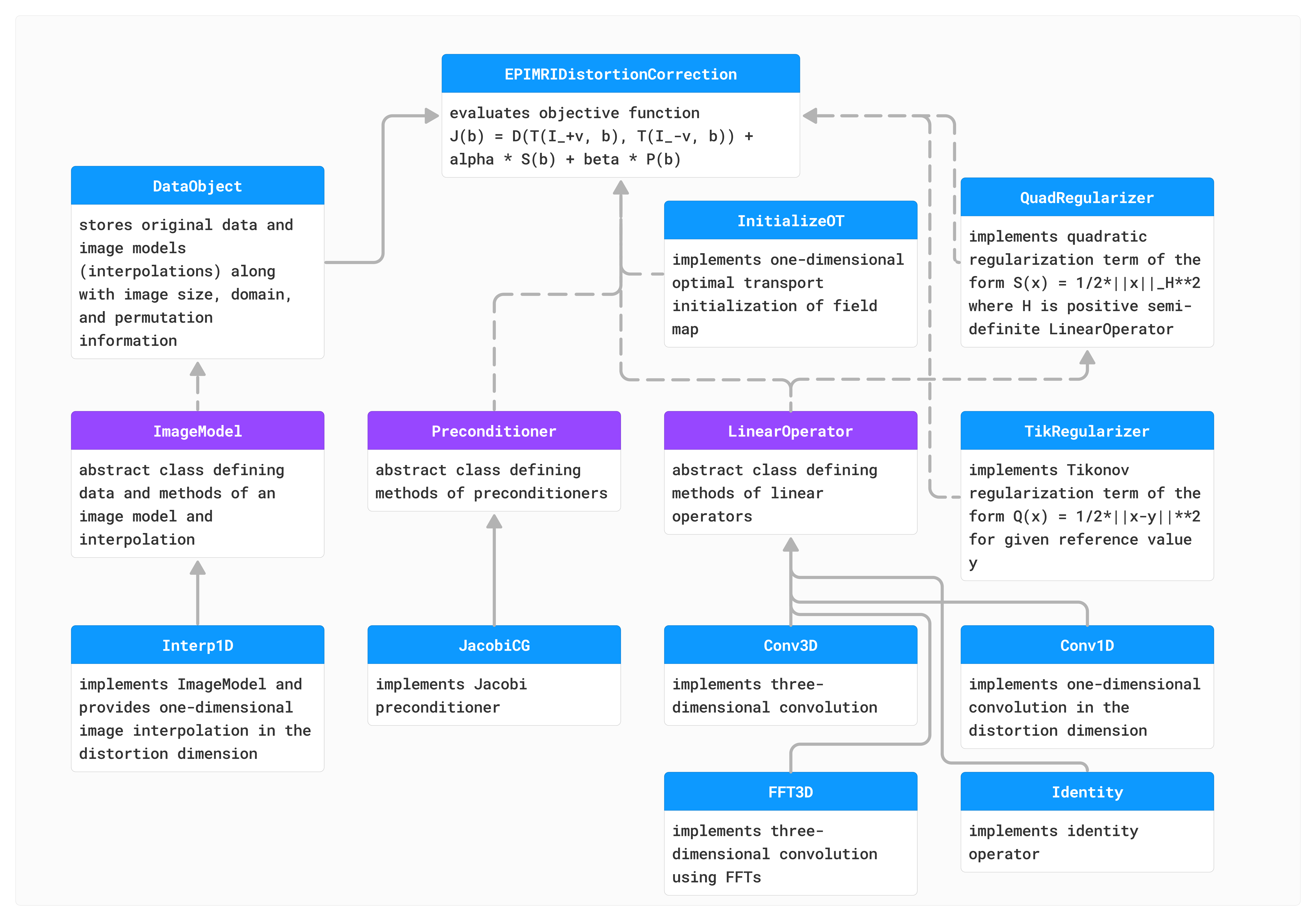}
        \caption{Class structure of PyHySCO loss function. The main class representing the loss function is \texttt{EPIMRIDistortionCorrection}. Purple classes are abstract, and blue classes are concrete. Solid arrows indicate inheritance. Dashed arrows indicate dependencies and class objects that are attributes.}
        \label{fig:UML_loss}
    \end{subfigure}
    \begin{subfigure}{0.9\textwidth}
            \includegraphics[width=\textwidth]{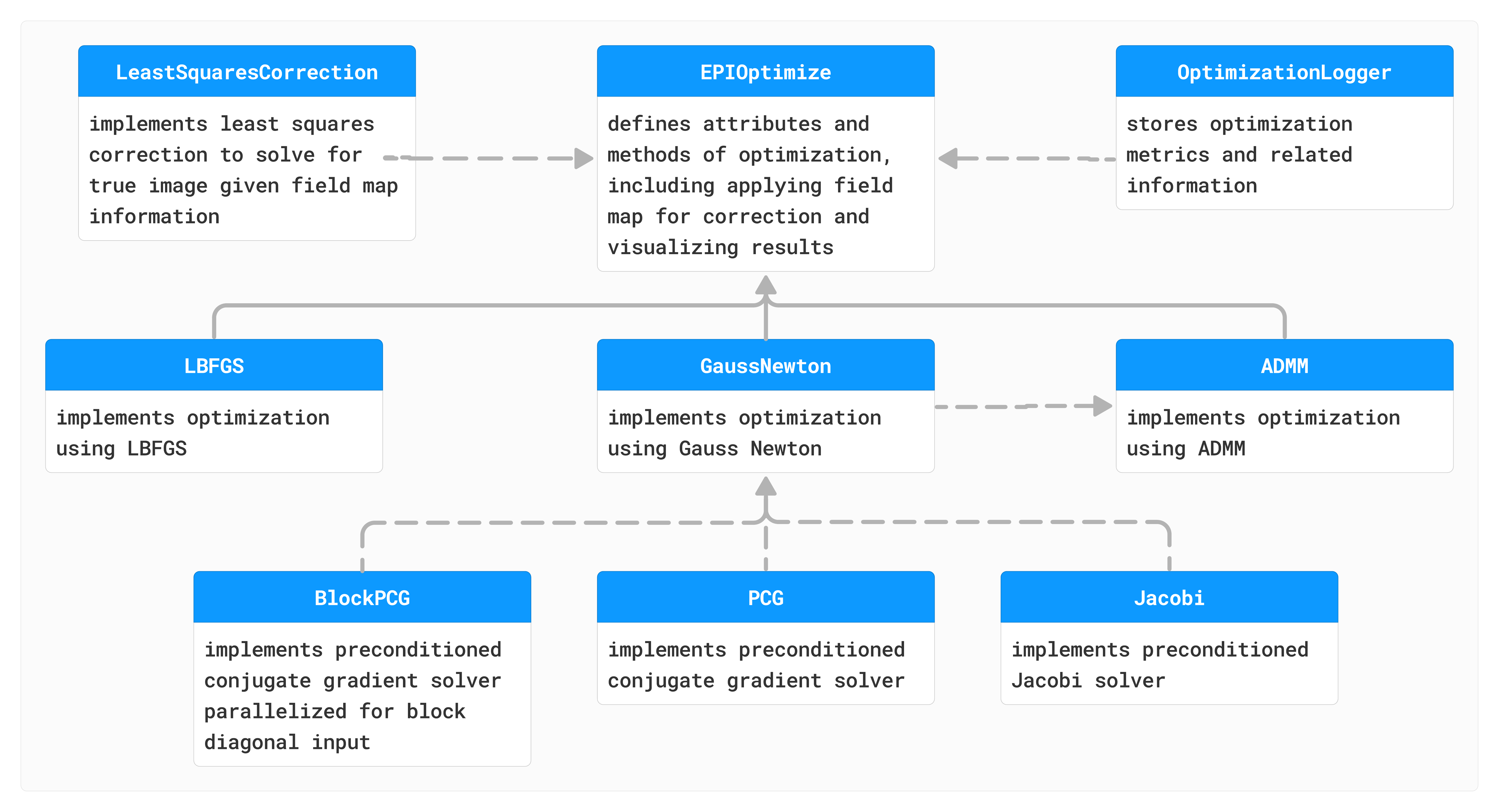}
    \caption{Class structure of PyHySCO optimization. The main class defining optimization is \texttt{EPIOptimize}. Solid arrows indicate inheritance. Dashed arrows indicate dependencies and class objects that are attributes.}
    \label{fig:UML_opt}
    \end{subfigure}
    \caption{UML diagram of PyHySCO showing the classes and relationships for the (\ref{fig:UML_loss}) loss function and (\ref{fig:UML_opt}) optimization. A \texttt{EPIMRIDistortionCorrection} object defining the loss function is an attribute of every \texttt{EPIOptimize} object defining the optimization scheme.}
    \label{fig:UML}
\end{figure}

We implemented PyHySCO in PyTorch \cite{pytorch} following the overall code structure visualized in the diagrams in Figures~\ref{fig:UML_loss} and \ref{fig:UML_opt} for the objective function and optimization, respectively. The main classes of PyHySCO are the loss function, implemented in \texttt{EPIMRIDistortionCorrection}, and the optimization, defined in \texttt{EPIOptimize}. The other classes and methods, described in detail below, implement the components of the loss function evaluation and optimization schemes.

\subsubsection{Data Storage and Image Model}
The input pair of images with opposite phase encoding directions are loaded and permuted such that the distortion dimension is the last, as this is where PyTorch expects the batch dimension for parallelizing operations. Information on the input images is stored in an object of type \texttt{DataObject}. This class stores information on the image size, domain, voxel size, how to permute the data back to the input order, and the \texttt{ImageModel} for each input image. The \texttt{ImageModel} abstract class defines the structure and required methods for an image model, including storing the original data and providing a method \texttt{eval} that returns the data interpolated on the given points. We provide the default implementation \texttt{Interp1D}, a piecewise linear one-dimensional interpolation parallelized in the last dimension. The \texttt{DataObject} for a given input pair is then stored in the \texttt{EPIMRIDistortionCorrection} object.

\subsubsection{Correction Model}
The mass-preserving correction model \eqref{eq:transform_model} is implemented in the method \texttt{mp\_transform}, a class method of  \texttt{EPIMRIDistortionCorrection}. The method takes as input an \texttt{ImageModel} and a field map. The geometric deformation is computed by using an averaging \texttt{LinearOperator} to compute the field map values in the cell centers and adding this to a cell-centered grid to obtain the deformed grid defined by this field map. Using the \texttt{ImageModel}, the image is interpolated on this deformed grid. The intensity modulation term is computed using a finite difference \texttt{LinearOperator}. The two terms are multiplied together element-wise before returning the corrected image. The default implementation of the \texttt{LinearOperator} objects for averaging and finite difference are given as one-dimensional convolutions, parallelized in the last dimension.

\subsubsection{Regularization Terms}
The intensity regularization term is computed within the \texttt{EPIMRIDistortionCorrection} class in the method \texttt{phi\_EPI} which computes the result of applying $\phi$ as defined in \eqref{eq:phi} element-wise to the result of applying the finite difference operator to the field map, as computed in the correction model. This function acts as a barrier term, ensuring that the derivative of the field map in the distortion dimension is in the range (-1, 1). 

The smoothness regularization term is implemented in a \texttt{QuadRegularizer} object, which defines the evaluation of a quadratic regularization term of the form of \eqref{eq:smoothness} using a positive semi-definite \texttt{LinearOperator} as $H$. By default, $H$ is a discretized negative Laplacian applied via a three-dimensional convolution.

In the ADMM optimizer, the regularizer structure differs to account for the splitting in~\eqref{eq:split}. The objective function for the $\mathbf{b}$ update in \eqref{eq:x_update} is computed in \texttt{EPIMRIDistortionCorrection} where the computation of $S_3$ is a one-dimensional Laplacian in the distortion dimension applied via a one-dimensional convolution. The proximal term is computed through a \texttt{TikRegularizer} object, a Tikhonov regularizer structure. The objective function for the $\mathbf{z}$ update in \eqref{eq:z_update} is a \texttt{QuadRegularizer} object where the \texttt{LinearOperator} $H$ is a two-dimensional Laplacian corresponding to $S_2$ and $S_3$. This operator is implemented in \texttt{FFT3D}, which defines an operator applying a convolution kernel diagonalized in Fourier space \cite{cooley1969convfft}. This implementation allows for easily inverting the kernel in solving for $\mathbf{z}$.

\subsubsection{Hessian and Preconditioning}
For the Gauss-Newton and ADMM optimizers, an approximate Hessian and preconditioner are additionally computed. Parts of the Hessian are computed in \texttt{EPIMRIDistortionCorrection} during objective function evaluation, and the Hessian can be applied through a matrix-vector product. Similarly, a \texttt{Preconditioner} can be computed during objective function evaluation and is accessible through a returned function applying the preconditioner to its input. By default, we provide a Jacobi preconditioner in the class \texttt{JacobiCG}.

\subsubsection{Initialization}
The \texttt{EPIMRIDistortionCorrection} class has a method \texttt{initialize}, returning an initial guess for the field map using some \texttt{InitializationMethod}. We provide an implementation of the proposed parallelized optimal transport-based initialization in \texttt{InitializeOT}. The implementation computes the one-dimensional transport maps in parallel using a linear spline interpolation. In practice, the parallelized initialization gives a highly non-smooth initial field map, so the method optionally applies a blurring operator using a 3-by-3-by-3 Gaussian kernel with a standard deviation of 1.0 to promote a smoother optimized field map. Applying the blur to the field map is implemented using the fast FFT convolution operator \texttt{FFT3D}.

\subsubsection{Optimization}
The minimization of the objective function defined in a \texttt{EPIMRIDistortionCorrection} object happens in a subclass of \texttt{EPIOptimize}, which takes the objective function object as input. During optimization, the \texttt{OptimizationLogger} class is used to track iteration history, saving it to a log file and optionally printing this information to standard output. PyHySCO includes implementations of the LBFGS, Gauss-Newton, and ADMM solvers described previously. Each of the classes \texttt{LBFGS}, \texttt{GaussNewton}, and \texttt{ADMM} provide a \texttt{run\_correction} method which minimizes the objective function using the indicated optimization scheme. The \texttt{LBFGS} implementation uses the explicitly computed derivative from \texttt{EPIMRIDistortionCorrection}. For LBFGS, we use the norm of the gradient reaching a given tolerance as stopping criteria, or the change in loss function or field map between iterations falling below a given tolerance. The \texttt{GaussNewton} implementation uses a conjugate gradient solver implemented in the class \texttt{PCG}. Our Gauss Newton implementation uses the same stopping criteria as LBFGS. The \texttt{ADMM} implementation solves the $\mathbf{b}$ update in \eqref{eq:x_update} using \texttt{GaussNewton} with a parallelized conjugate gradient solver in \texttt{BlockPCG}. The $\mathbf{z}$ update in \eqref{eq:z_update} is solved directly through the inverse method \texttt{inv} of the operator used to define the \texttt{QuadRegularizer} for this term, efficiently implemented using FFTs in \texttt{FFT3D}. As stopping criteria, the ADMM iterations will terminate if the change in all of $\mathbf{b}$, $\mathbf{z}$, and $\mathbf{u}$ from the previous iteration falls below a given tolerance.

\subsubsection{Image Correction}
The optimal field map, stored as \texttt{Bc} in the \texttt{EPIOptimize} object after \texttt{run\_correction} is completed, can be used to produce a corrected image or pair of images. The \texttt{apply\_correction} method of \texttt{EPIOptimize} implements both a Jacobian modulation correction and a least squares correction. The Jacobian modulation correction is based on the model of \cite{ChangFitzpatrick1992} as implemented in the \texttt{mp\_transform} method of \texttt{EPIMRIDistortionCorrection}. This correction method computes and saves two corrected images, one for each input image.

The field map can also be used in a least squares correction similar to the correction in \cite{AnderssonEtAl2003}, implemented in \texttt{LeastSquaresCorrection}. In this correction, the estimated field map determines a push-forward matrix that transforms the true image to the distorted image given as input. This gives rise to a least squares problem for the true image given the input images and push forward matrix.

\subsection{PyHySCO Usage and Workflow}
\label{sub:usage}
\begin{figure}
    \centering
    \begin{subfigure}{0.9\textwidth}
        \includegraphics[width=\textwidth]{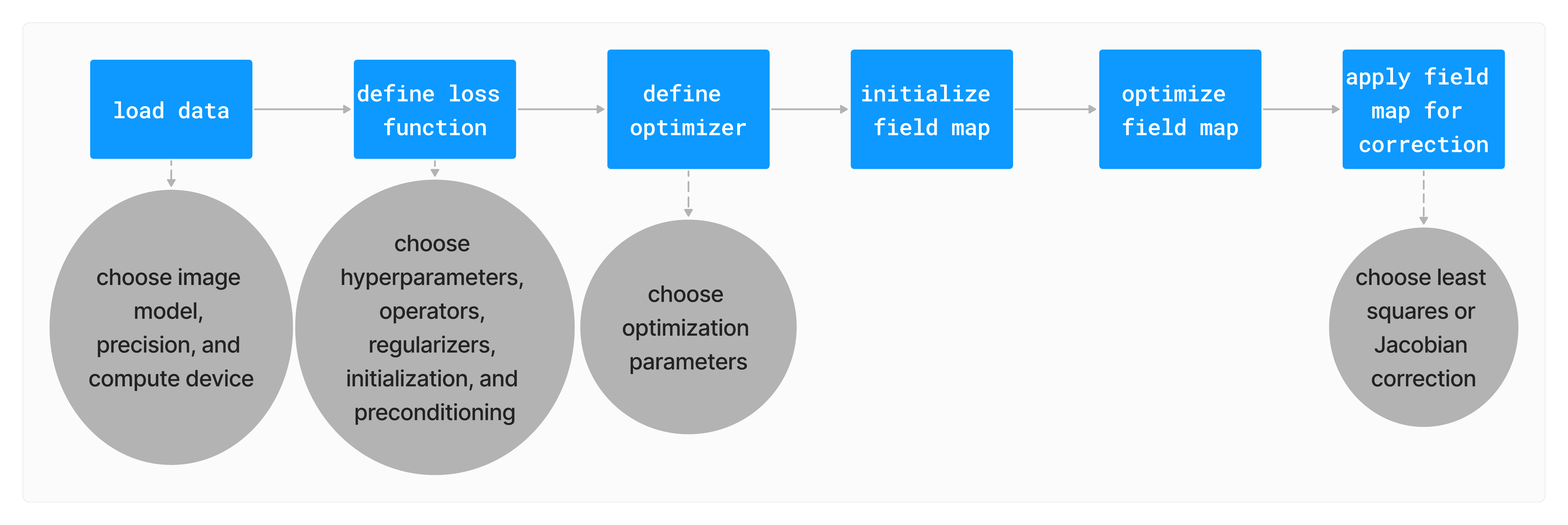}
        \caption{The workflow of the PyHySCO toolbox from setup through optimization and distortion correction.}
        \label{fig:flowchart}
    \end{subfigure}
    \begin{subfigure}{0.9\textwidth}
            \includegraphics[width=\textwidth]{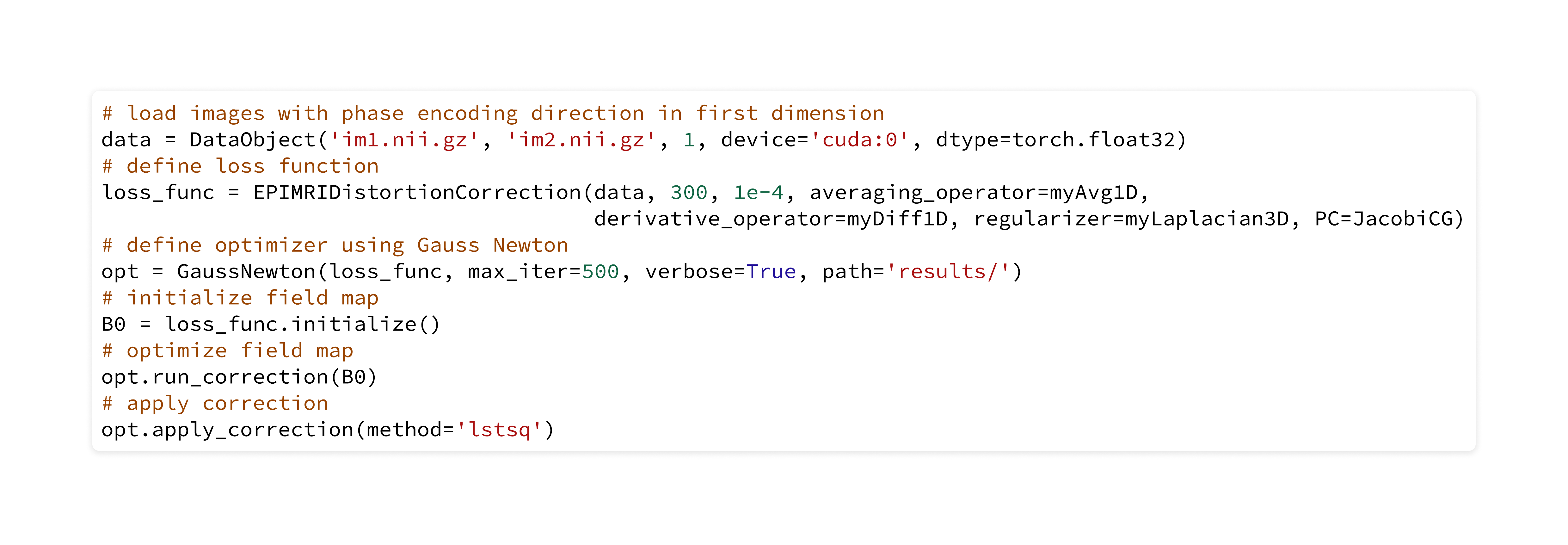}
    \caption{An example using the PyHySCO toolbox from a Python script.}
    \label{fig:script}
    \end{subfigure}
    \begin{subfigure}{0.9\textwidth}
            \includegraphics[width=\textwidth]{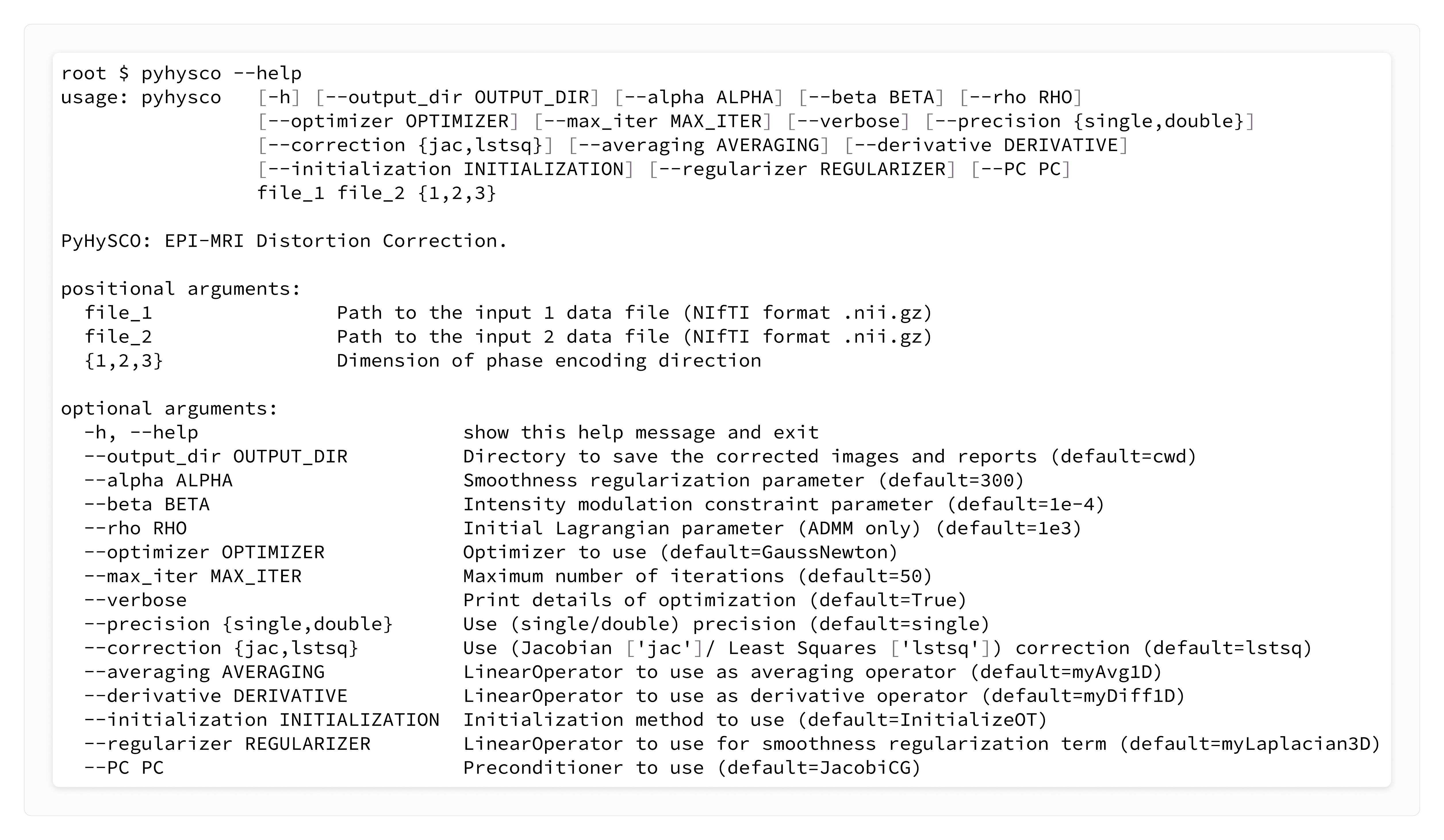}
    \caption{The help message for the PyHySCO command line interface. This interface allows use of PyHySCO as part of existing MRI post-processing pipelines.}
    \label{fig:commandline}
    \end{subfigure}
    \caption{The usage and workflow of PyHySCO.}
    \label{fig:walkthrough}
\end{figure}

The workflow of PyHySCO is illustrated in Figure \ref{fig:flowchart} alongside examples of using PyHySCO in a Python script (Figure \ref{fig:script}) and through the command line (Figure \ref{fig:commandline}). Running PyHySCO from a user-defined Python script allows more specific control of the inputs and outputs from PyHySCO methods. The command line interface allows the user to pass configuration options directly from the command-line, which enables our EPI distortion correction tool to be easily used as a part of existing command-line based MRI post-processing pipelines such as the FSL toolbox \cite{smith2004FSL}. Executing PyHySCO requires the user to provide, at a minimum, the file paths for the input pair of images with opposite phase encoding directions and which dimension (1, 2, or 3) is aligned with the phase encoding direction. The modularity of PyHySCO additionally allows for configuring options such as the scalar hyperparameters in \eqref{eq:opt_disc}, implementation of operators, regularizers, and interpolation, optimizer and associated optimization parameters, and image correction method.

Regardless of execution through a script or the command line, PyHySCO stores the input images in a \texttt{DataObject} object, the loss function in a \texttt{EPIMRIDistortionCorrection} object, and the optimizer in an object of a subclass of \texttt{EPIOptimize}. The field map is initialized from the method \texttt{initialize} in \texttt{EPIMRIDistortionCorrection}, and the field map is optimized by calling the method \texttt{run\_correction} in the optimizer object. Finally, the method \texttt{apply\_correction} in \texttt{EPIOptimize} applies the field map to correct the input images and saves the result to one or more NIFTI file(s).

\section{Results}\label{sec:3results}
We demonstrate PyHySCO's effectiveness through extensive experiments using real and simulated data from the Human Connectome Project \cite{HCPdata} and validate the novel initialization scheme and implementation of optimization algorithms. 
Section~\ref{subs:validation} describes the datasets and Section~\ref{subs:ethics} contains our ethics statement.
Section~\ref{subs:metrics} introduces our evaluation metrics.
The experiments in~\ref{subs:optimizers} compare the performance of the three optimization algorithms implemented in PyHySCO on CPU and GPU hardware.
In Section~\ref{subs:precision}, we empirically show that PyHySCO can be effective in single precision, which accelerates computation on modern GPU hardware.
The results in Section~\ref{subs:HyscoTopup}, suggest that PyHySCO can achieve correction quality at a considerably shorter runtime compared to  HySCO and TOPUP \cite{RuthottoEtAl2013hysco, AnderssonEtAl2003}.

\subsection{Validation Datasets}
\label{subs:validation}
The data used in the following experiments is from the Human Connectome Project \cite{HCPdata}. We validate our methods and tool on 3T and 7T diffusion-weighted imaging data from the HCP 1200 Subjects Release, with 20 subjects randomly chosen for each field strength. Table \ref{table:HCPdata} provides details of the datasets.

We also evaluate our methods on simulated data. This data only contains susceptibility artifact distortions, so it shows how our tool performs without the influence of other factors, e.g., patient movement between scans. To simulate the distortions, we use a pair of magnitude and phase images for a subject in HCP and generate the field map using FSL's FLIRT and PRELUDE tools \cite{smith2004FSL}. Considering the physical model of \cite{ChangFitzpatrick1992}, the field map $b$ can be used to define the push-forward matrices that give how the intensity value at $x$ is pushed forward to $x+b(x)$ in the distortion direction $+v$ as well as the opposite direction $-v$. Applying the push-forward matrices to a T2-weighted image for the subject, we generate a pair of distorted images. For the simulated data, we then have a reference value for the field map and an undistorted, true image.

\subsection{Ethics Statement}
\label{subs:ethics}
No new data were collected specifically for this paper. The Human Connectome Project, through the WU-Minn HCP Consortium, obtained written informed consent from all participants in the data collection study \cite{HCPdata}.

\begin{table}[t]
\begin{center}
\begin{tabular}{||c || c | c | c | c ||} 
 \hline
 Dataset & No. of Subjects & Image Size & Resolution & PE directions \\ [0.5ex] 
 \hline\hline
 3T & 20 & 168 $\times$ 144 $\times$ 111 & 1.25 $\times$ 1.25 $\times$ 1.25 $mm^3$ &  LR/ RL \\
 \hline
 7T & 20 & 200 $\times$ 200 $\times$ 132 & 1.05 $\times$ 1.05 $\times$ 1.05 $mm^3$ &  AP/ PA \\
 \hline
 Simulated & 20 & 320 $\times$ 320 $\times$ 256 & 0.7 $\times$ 0.7 $\times$ 0.7 $mm^3$ & AP/ PA \\
 \hline
\end{tabular}
 \caption{Details of data used in validation. LR/RL is left-to-right and right-to-left phase encoding, and AP/PA is anterior-to-posterior and posterior-to-anterior phase encoding. Further details of acquisition parameters are in \cite{HCPdata}.}
 \label{table:HCPdata} 
\end{center}
\end{table}

\subsection{Metrics and Comparison Methods}
\label{subs:metrics}
The quality of correction results is measured using the relative improvement of the distance between a pair of corrected images. Particularly, we calculate the sum-of-squares distance (SSD) of the corrected image pair relative to the SSD of the input pair. This metric is a useful surrogate for the correctness of the field map in the absence of a ground truth \cite{graham2017quantitative}. Additionally, we take the value of the smoothness regularization term $S(\mathbf{b})$ as a measure of how smooth the resulting field map is, with lower values being better.

We report the runtime in seconds of PyHySCO. The runtime is measured as the wall clock time using the Linux \texttt{time} command when calling the correction method from the command line. This time, therefore, includes the time taken to load and save the image data. In some cases, we also report the optimization time only, without loading and saving data, as measured by Python's \texttt{time} module.

We compare the runtime, relative improvement, and resulting images after correction using PyHySCO against those given by TOPUP \cite{AnderssonEtAl2003} as implemented in FSL \cite{smith2004FSL} using the default configuration\footnote{The default TOPUP configuration performs upsampling requiring the dimensions to be a multiple of 2. The configuration for TOPUP with images of the 3T data set does not perform upsampling due to the odd number of slices in the image volumes.}, and HySCO \cite{RuthottoEtAl2013hysco} as implemented in the ACID toolbox for SPM using the default parameters. HySCO is also based on the optimization problem \eqref{eq:opt_disc}, while TOPUP uses a slightly different objective function. This makes it difficult to compute smoothness and loss function values for TOPUP.

\subsection{Validity of Optimal Transport Initialization}
\label{subs:OTinit}
\begin{figure}[t]
    \centering
    \input{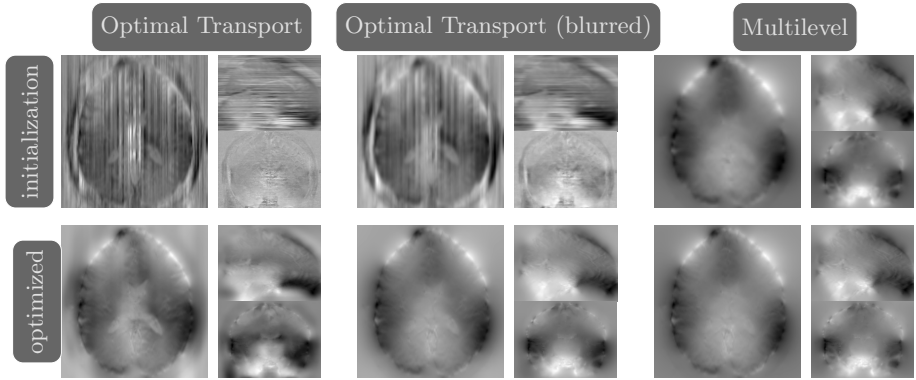}
    \caption{Example field maps (Subject ID 826353) at initialization (top row) and after optimization with Gauss-Newton (bottom row). The first column uses the proposed optimal transport initialization scheme. The middle column uses the same scheme with an additional Gaussian blur to promote smoothness. The right column uses the coarse-to-fine multilevel initialization scheme from HySCO with five levels, and the final field map is optimized at the original image resolution. The multilevel initialized field map is smooth by construction and further optimized to improve the relative image distance at the full resolution. The optimal transport initialization accurately corrects the distortions but is not smooth in the non-distortion dimensions unless blurred with a Gaussian. After the fine-level optimization all field maps are visually similar. }
    \label{fig:resultsInit}
\end{figure}

We compare the results of PyHySCO using our optimal transport initialization to those of the multi-level initialization used in HySCO  \cite{RuthottoEtAl2013hysco} both at initialization and after optimization with Gauss-Newton. The multi-level optimization of HySCO solves the optimization problem on a coarse grid and uses the result as the initialization of optimization on a finer grid, continuing until reaching the original image resolution; this follows the guidelines of~\cite[Chapter 9.4]{Modersitzki2009}. In our experiments, we use five levels in the initialization. The multi-level initialization gives a field map that is smooth by construction and improves the distance reduction as the grid becomes more fine. The field map from the PyHySCO optimal transport initialization drastically lowers the relative error between the input images, a relative improvement of over 96\% on real data and 94\% on simulated data. However, the parallelized one-dimensional computations lead to a lack of smoothness in the resulting field map. The smoothness can be improved by applying a Gaussian blur to the field map from the optimal transport initialization. This field map is smoother after initialization and gives a smoother field map after optimization. These results are comparable in relative error and smoothness to the field map optimized from the multilevel initialization of HySCO. Our one-dimensional parallelized optimal transport, even with the additional Gaussian blur, is much faster to compute than the multilevel initial field map given the ability to parallelize computations. PyHySCO initialization on a GPU with the additional blur takes less than 1 second on real data and about 3 seconds on simulated data. In comparison, the multi-level initialization on a CPU takes 30 to 40 seconds on real data and over 2 minutes on simulated data. The mean and standard deviation relative improvement, smoothness value, loss function value, and runtime are reported in Table~\ref{table:initialization} across all datasets. Examples of these field maps before and after optimization are shown in Figure~\ref{fig:resultsInit}.

\begin{table}[t]
\begin{center}
\resizebox{\columnwidth}{!}{%
\begin{tabular}{||c|c||c|c||c|c||c|c||}
\hline
\multicolumn{2}{||c||}{} & \multicolumn{2}{c||}{Optimal Transport} & \multicolumn{2}{c||}{Optimal Transport (blur)} & \multicolumn{2}{c||}{Multilevel}\\
\cline{3-8}
\multicolumn{2}{||c||}{} & initial & after opt & initial & after opt & initial & after opt \\
\hline\hline
\multirow{10}{*}{3T} & \multirow{2}{*}{Runtime (s)} & $5.78$ & $11.43$ & $6.31$ & $15.36$ & $41.69$ & $55.34$ \\
 & & $\pm 1.26$ & $\pm 1.46$ & $\pm 0.60$ & $\pm 3.90$ & $\pm 1.71$ & $\pm 2.84$ \\ \cline{2-8}
 & \multirow{2}{*}{Opt. Time (s)} &  $0.27$ & $4.34$ & $0.28$ & $6.78$ & $42.43$ & $48.65$ \\
 & & $\pm 0.01$ & $\pm 0.67$ & $\pm 0.02$ & $\pm 0.67$ & $\pm 4.04$ & $\pm 3.95$ \\ \cline{2-8}
 & Relative & $96.44$ & $83.90$ & $79.71$ & $82.75$ & $67.04$ & $81.96$ \\
 & Improvement & $\pm 1.13$ & $\pm 3.43$ & $\pm 3.43$ & $\pm 3.49$ & $\pm 5.15$ & $\pm 3.51$ \\ \cline{2-8}
 & \multirow{2}{*}{Loss Value} & $1.05e09$ & $2.84e07$ & $1.76e08$ & $2.56e07$ & $4.82e07$ & $2.51e07$ \\
 & & $\pm 2.66e08$ & $\pm 7.49e06$ & $\pm 5.20e07$ & $\pm 7.66e06$ & $\pm 1.70e07$ & $\pm 7.54e06$ \\ \cline{2-8}
 & Smoothness & $3.50e06$ & $5.08e04$ & $5.28e05$ & $3.85e04$ & $6.89e04$ & $3.47e04$ \\
 & Reg. Value & $\pm 8.81e05$ & $\pm 1.23e04$ & $\pm 1.54e05$ & $\pm 1.21e04$ & $\pm 2.98e04$ & $\pm 1.08e04$ \\ \cline{2-8}
\hline\hline
\multirow{10}{*}{7T} & \multirow{2}{*}{Runtime (s)} & $7.61$ & $13.55$ & $8.32$ & $19.72$ & $58.73$ & $77.79$ \\
 & & $\pm 1.99$ & $\pm 2.04$ & $\pm 2.79$ & $\pm 2.91$ & $\pm 6.24$ & $\pm 5.72$ \\ \cline{2-8}
  & \multirow{2}{*}{Opt. Time (s)} & $0.61$ & $5.09$ & $0.63$ & $10.16$ & $30.38$ & $40.50$ \\
 & & $\pm 0.02$ & $\pm 1.23$ & $\pm 0.02$ & $\pm 0.85$ & $\pm 2.63$ & $\pm 3.87$ \\ \cline{2-8}
 & Relative & $96.53$ & $86.01$ & $75.09$ & $85.76$ & $69.12$ & $85.42$ \\
 & Improvement & $\pm 1.47$ & $\pm 5.15$ & $\pm 3.97$ & $\pm 5.10$ & $\pm 8.28$ & $\pm 5.08$ \\ \cline{2-8}
 & \multirow{2}{*}{Loss Value} & $3.48e09$ & $5.28e07$ & $4.50e08$ & $4.14e07$ & $7.77e07$ & $4.02e07$ \\
 & & $\pm 1.15e09$ & $\pm 2.01e07$ & $\pm 2.47e08$ & $\pm 1.95e07$ & $\pm 3.01e07$ & $\pm 1.82e07$ \\ \cline{2-8}
 & Smoothness & $1.16e07$ & $9.52e04$ & $1.36e06$ & $5.63e04$ & $8.21e04$ & $5.03e04$ \\
 & Reg. Value & $\pm 3.83e06$ & $\pm 2.64e04$ & $\pm 7.74e05$ & $\pm 1.91e04$ & $\pm 4.07e04$ & $\pm 1.48e04$ \\ \cline{2-8}
\hline\hline
\multirow{10}{*}{Simulated} & \multirow{2}{*}{Runtime (s)} & $10.62$ & $80.29$ & $16.59$ & $106.47$ & $173.20$ & $47.98$ \\
 & & $\pm 0.57$ & $\pm 9.96$ & $\pm 0.64$ & $\pm 11.21$ & $\pm 27.06$ & $\pm 8.38$ \\ \cline{2-8}
  & \multirow{2}{*}{Opt. Time (s)} & $3.51$ & $64.45$ & $3.61$ & $89.17$ & $125.35$ & $157.95$ \\
 & & $\pm 0.03$ & $\pm 10.02$ & $\pm 0.15$ & $\pm 11.48$ & $\pm 24.88$ & $\pm 28.79$  \\ \cline{2-8}
 & Relative & $94.64$ & $76.82$ & $75.34$ & $76.27$ & $55.01$ & $73.63$ \\
 & Improvement & $\pm 1.26$ & $\pm 5.09$ & $\pm 3.44$ & $\pm 5.18$ & $\pm 5.66$ & $\pm 5.39$ \\ \cline{2-8}
 & \multirow{2}{*}{Loss Value} & $5.10e08$ & $6.31e07$ & $2.11e08$ & $6.07e07$ & $8.17e07$ & $5.83e07$ \\
 & & $\pm 9.51e07$ & $\pm 1.46e07$ & $\pm 4.30e07$ & $\pm 1.39e07$ & $\pm 2.08e07$ & $\pm 1.33e07$ \\ \cline{2-8}
 & Smoothness & $1.67e06$ & $1.06e05$ & $5.84e05$ & $9.53e04$ & $6.18e04$ & $7.50e04$ \\
 & Reg. Value & $\pm 3.14e05$ & $\pm 2.94e04$ & $\pm 1.24e05$ & $\pm 2.71e04$ & $\pm 1.70e04$ & $\pm 2.20e04$ \\ \cline{2-8}
\hline\hline
\end{tabular}
}%
\caption{Validation of the optimal transport initialization. We compare the runtime, relative improvement, smoothness value, and loss function value at initialization and after optimization with Gauss Newton for the proposed parallelized optimal transport initializaztion, the proposed initialization with an additional Gaussian blur, and the multilevel initialization used in HySCO \cite{RuthottoEtAl2013hysco}. For each metric we report the mean and standard deviation in the 3T, 7T, and simulated datasets. The multilevel initialization is timed on CPU in Matlab, and the optimal transport initializations and all optimizations are timed on GPU in Python. The optimal transport based initializations provide a comparable quality while decreasing runtime compared to the multilevel initialization, and the optimal transport with Gaussian blur promotes a more smooth field map.}
\label{table:initialization}
\end{center}
\end{table}

\subsection{Comparison of PyHySCO Optimizers on GPU and CPU}
\label{subs:optimizers}
We compare the results of PyHySCO using GN-PCG, ADMM, and LBFGS on both GPU and CPU architectures. Table~\ref{table:optimizerscpugpu} shows the runtimes and correction quality of each optimizer on CPU and GPU. All optimizers achieve a similar correction quality with respect to relative improvement of image distance, loss value, and smoothness regularizer value, but GN-PCG has faster runtime on both CPU and GPU. On real data, GN-PCG took 10-13 seconds on average on GPU and 27-31 seconds on average on CPU, while ADMM took 11-15 seconds on GPU and 98-158 seconds on CPU, and LBFGS took 23-36 seconds on GPU and 104-141 seconds on CPU. Table~\ref{table:optimizersmetrics} shows optimization metrics, including the number of iterations, stopping criteria, number of function evaluations, number of Hessian evaluations, and number of inner iterations if applicable. Consistent with its faster runtime, optimization with GN-PCG achieves a similar loss value with less computation as measured by function and Hessian evaluations. Figures~\ref{fig:results3T}, \ref{fig:results7T}, and~\ref{fig:resultssim} show the field map and corrected images for each optimizer for one example subject from each dataset. The field maps and corrected images are visually similar across optimizers.

\begin{table}[t]
\begin{center}
\begin{tabular}{||c|c||c|c||c|c||c|c||}
\hline
\multicolumn{2}{||c||}{} & \multicolumn{2}{c||}{LBFGS} & \multicolumn{2}{c||}{GN-PCG} & \multicolumn{2}{c||}{ADMM}\\
\cline{3-8}
\multicolumn{2}{||c||}{} & CPU & GPU & CPU & GPU & CPU & GPU \\
\hline\hline
\multirow{10}{*}{3T} & \multirow{2}{*}{Runtime (s)} & $104.45$ & $23.13$ & $27.37$ & $10.37$ & $98.54$ & $11.58$ \\
 & & $\pm 70.74$ & $\pm 4.61$ & $\pm 4.53$ & $\pm 0.87$ &  $\pm 30.15$ & $\pm 2.23$ \\ \cline{2-8}
 & \multirow{2}{*}{Opt. Time (s)} & $100.28$ & $16.70$ & $23.13$ & $4.38$ & $94.53$ & $5.63$ \\
 & & $\pm 70.82$ & $\pm 4.49$ & $\pm 4.53$ & $\pm 0.68$ & $\pm 30.20$ & $\pm 2.15$ \\ \cline{2-8}
 & Relative & $81.47$ & $82.32$ & $82.74$ & $82.74$ &  $82.76$ & $82.77$ \\
 & Improvement & $\pm 3.71$ & $\pm 3.40$ & $\pm 3.50$ & $\pm 3.50$ & $\pm 3.31$ & $\pm 3.30$ \\ \cline{2-8}
 & \multirow{2}{*}{Loss Value} & $7.90e07$ & $2.56e07$ & $2.56e07$ & $2.56e07$ & $3.09e07$ & $3.10e07$ \\
 & & $\pm 7.99e07$ & $\pm 7.72e06$ & $\pm 7.69e06$ & $\pm 7.69e06$ & $\pm 8.51e96$ & $\pm 8.56e06$ \\ \cline{2-8}
 & Smoothness & $2.13e05$ & $3.72e04$ & $3.85e04$ & $3.85e04$ & $5.62e04$ & $5.65e04$ \\
 & Reg. Value & $\pm 2.56e05$ & $\pm 1.18e04$ & $\pm 1.21e04$ & $\pm 1.21e04$  & $\pm 1.65e04$ & $\pm 1.71e04$ \\ \cline{2-8}
\hline\hline
\multirow{10}{*}{7T} & \multirow{2}{*}{Runtime (s)} & $141.44$ & $36.23$ & $31.71$ & $13.62$ & $158.64$ & $15.25$ \\
 & & $\pm 117.38$ & $\pm 7.76$ & $\pm 3.18$ & $\pm 2.38$ & $\pm 46.99$ & $\pm 3.15$ \\ \cline{2-8}
  & \multirow{2}{*}{Opt. Time (s)} & $135.72$ & $29.23$ & $26.84$ & $6.57$ & $152.69$ & $8.34$ \\
 & & $\pm 116.29$ & $\pm 7.88$ & $\pm 3.15$ & $\pm 2.30$ & $\pm 46.64$ & $\pm 2.91$ \\ \cline{2-8}
 & Relative & $80.75$ & $85.74$ & $85.76$ & $85.76$ & $85.87$ & $85.85$ \\
 & Improvement & $\pm 6.91$ & $\pm 4.99$ & $\pm 5.10$ & $\pm 5.10$ & $\pm 4.99$ & $\pm 4.99$ \\ \cline{2-8}
 & \multirow{2}{*}{Loss Value} & $2.25e08$ & $4.25e07$ & $4.14e07$ & $4.14e07$ & $4.43e07$ & $4.43e07$ \\
 & & $\pm 2.22e08$ & $\pm 2.00e07$ & $\pm 1.95e07$ & $\pm 1.95e07$ & $\pm 1.99e07$ & $\pm 1.95e07$ \\ \cline{2-8}
 & Smoothness & $6.38e05$ & $6.00e04$ & $5.63e04$ & $5.63e04$ & $6.68e04$ & $6.66e04$ \\
 & Reg. Value & $\pm 7.01e05$ & $\pm 2.18e04$ & $\pm 1.91e04$ & $\pm 1.91e04$ & $\pm 2.18e04$ & $\pm3.68e04$ \\ \cline{2-8}
\hline\hline
\multirow{10}{*}{Sim.} & \multirow{2}{*}{Runtime (s)} & $6344.93$ & $143.77$ & $1094.96$ & $55.26$ & $7687.28$ & $52.72$ \\
 & & $\pm 649.21$ & $\pm 6.47$ & $\pm 135.20$ & $\pm 3.86$ & $\pm 4596.31$ & $\pm 18.01$ \\ \cline{2-8}
  & \multirow{2}{*}{Opt. Time (s)} & $6320.43$ & $125.95$ & $1070.65$ & $37.60$ & $7662.55$ & $35.15$ \\
 & & $\pm 649.01$ & $\pm 6.40$ & $\pm 135.69$ & $\pm 4.54$ & $\pm 4596.38$ & $\pm 17.92$ \\ \cline{2-8}
 & Relative & $75.45$ & $75.44$ & $76.28$ & $76.28$ & $74.93$ & $75.00$ \\
 & Improvement & $\pm 5.40$ & $\pm 5.35$ & $\pm 5.19$ & $\pm 5.18$ & $\pm 5.59$ & $\pm 5.34$ \\ \cline{2-8}
 & \multirow{2}{*}{Loss Value} & $6.03e07$ & $6.00e07$ & $6.08e07$ & $6.08e07$ & $6.08e07$ & $6.12e07$ \\
 & & $\pm 1.44e07$ & $\pm 1.41e07$ & $\pm 1.40e07$ & $\pm 1.40e07$ & $\pm 1.40e07$ & $\pm 1.43e07$ \\ \cline{2-8}
 & Smoothness & $9.06xe04$ & $8.94e04$ & $9.56e04$ & $9.56e04$ & $8.97e04$ & $9.12e04$ \\
 & Reg. Value & $\pm 2.94e04$ & $\pm 2.74e04$ & $\pm 2.74e04$ & $\pm 2.72e04$ & $\pm 2.79e04$ & $\pm 2.77e04$ \\ \cline{2-8}
\hline\hline
\end{tabular}
\caption{The speed and quality of optimization in PyHySCO on GPU and CPU with LBFGS, Gauss Newton, and ADMM. We report for each dataset and optimizer the mean and standard deviation total runtime (including loading and saving data), optimization time, improvement in distance between corrected images relative to input image, loss value, and smoothness regularizer value. Gauss Newton achieves a similar correction quality in less time than LBFGS or ADMM on both CPU and GPU.}
\label{table:optimizerscpugpu}
\end{center}
\end{table}
\begin{table}[t]
\begin{center}
\begin{tabular}{||c|c||c||c||c||}
\hline
\multicolumn{2}{||c||}{} & \multicolumn{1}{c||}{LBFGS} & \multicolumn{1}{c||}{GN-PCG} & \multicolumn{1}{c||}{ADMM}\\
\cline{3-5}
\hline\hline
\multirow{12}{*}{3T} & \multirow{2}{*}{Iterations} & $455.30$ & $8.400$ & $36.05$  \\
 & & $\pm 52.80$ & $\pm 0.92$ & $\pm 10.37$  \\ \cline{2-5}
  & Stopping Criteria & \multirow{2}{*}{9/3/0/8} & \multirow{2}{*}{0/20/0/0} & \multirow{2}{*}{0/0/20/0}  \\
 & (grad/loss/field map/max iter) &  & &  \\ \cline{2-5}
 & \multirow{2}{*}{Func. Evals} & $463.30$ & $9.40$ & $37.05$ \\
 & &$\pm 54.12$ & $\pm 0.92$ & $\pm 10.37$  \\ \cline{2-5}
 & \multirow{2}{*}{Hessian Evals} & \multirow{2}{*}{N/A} & $92.40$ & $437.50$  \\
 & & & $\pm 10.08$ & $\pm 140.30$  \\ \cline{2-5}
 & \multirow{2}{*}{Inner Iterations} &  \multirow{2}{*}{N/A} & $10.0000$ & $11.0269$ \\
 & & & $\pm 0.00$ & $\pm 1.02$   \\ \cline{2-5}
 & \multirow{2}{*}{Loss Value} & $2.56e07$ & $2.56e07$ & $3.10e07$ \\
 & & $\pm 7.72e06$ & $\pm 7.69e06$ & $\pm 8.56e06$ \\ \cline{2-5}
\hline\hline

\multirow{12}{*}{7T} & \multirow{2}{*}{Iterations} & $405.00$ & $7.50$ & $56.75$  \\
 & & $\pm 65.61$ & $\pm 0.87$ & $\pm 17.01$  \\ \cline{2-5}
  & Stopping Criteria & \multirow{2}{*}{14/3/0/3} & \multirow{2}{*}{0/20/0/0} & \multirow{2}{*}{0/0/20/0}  \\
 & (grad/loss/field map) &  & &  \\ \cline{2-5}
 & \multirow{2}{*}{Func. Evals} & $415.35$ & $8.50$ & $57.75$ \\
 & &$\pm 68.00$ & $\pm 0.87$ & $\pm 17.01$  \\ \cline{2-5}
 & \multirow{2}{*}{Hessian Evals} & \multirow{2}{*}{N/A} & $82.25$ & $339.05$  \\
 & & &  $\pm 9.15$ & $\pm 101.55$  \\ \cline{2-5}
 & \multirow{2}{*}{Inner Iterations} &  \multirow{2}{*}{N/A} & $9.9722$ & $4.9771$ \\
 & & & $\pm 0.12$ & $\pm 0.08$  \\ \cline{2-5}
 & \multirow{2}{*}{Loss Value}  & $4.25e07$ & $4.14e07$ & $4.43e07$ \\
 & & $\pm 2.00e07$ & $\pm 1.95e07$ & $\pm 1.95e07$ \\ \cline{2-5}
\hline\hline

\multirow{12}{*}{Simulated} & \multirow{2}{*}{Iterations} & $497.65$ & $20.05$ & $109.35$  \\
 & & $\pm 5.88$ & $\pm 1.83$ & $\pm 64.52$  \\ \cline{2-5}
  & Stopping Criteria & \multirow{2}{*}{1/0/0/19} & \multirow{2}{*}{0/18/2/0} & \multirow{2}{*}{0/0/20/0}  \\
 & (grad/loss/field map) &  & &  \\ \cline{2-5}
 & \multirow{2}{*}{Func. Evals} & $532.35$ & $21.05$ & $110.35$ \\
 & &$\pm 28.27$ & $\pm 1.83$ & $\pm 64.52$  \\ \cline{2-5}
 & \multirow{2}{*}{Hessian Evals} & \multirow{2}{*}{N/A} & $220.55$ & $1872.15$  \\
 & & & $\pm 20.13$ & $\pm 1417.11$  \\ \cline{2-5}
 & \multirow{2}{*}{Inner Iterations} &  \multirow{2}{*}{N/A} & $10.0000$ & $15.1681$ \\
 & & & $\pm 0.00$ & $\pm 3.69$ \\ \cline{2-5}
 & \multirow{2}{*}{Loss Value} & $6.00e07$ & $6.08e07$ & $6.12e07$ \\
 & & $\pm 1.41e07$ & $\pm 1.40e07$ & $\pm 1.43e07$ \\ \cline{2-5}
\hline\hline

\end{tabular}
\caption{Details of optimization for PyHySCO optimizers LBFGS, Gauss Newton, and ADMM. For each dataset we report the average and standard deviation number of iterations, count of stopping criteria used (gradient tolerance/ loss function change tolerance/ field map change tolerance/ maximum iterations), average and standard deviation number of function evaluations, average and standard deviation number of Hessian evaluations, average and standard deviation number of inner iterations, and average and standard deviation loss value. Gauss Newton achieves a similar quality of correction with less computation than LBFGS or ADMM.}
\label{table:optimizersmetrics}
\end{center}
\end{table}

\subsection{Single Precision vs Double Precision on GPU and CPU}
\label{subs:precision}
We show the validity of PyHySCO using optimal transport initialization and GN-PCG in both double precision (64 bit) and single precision (32 bit) arithmetic on three different GPU architectures and a CPU architecture. Since GPU architectures are optimized for the speed of lower precision calculations, we see a significant speedup when using single precision instead of double precision. Calculations in single precision, however, have the risk of lower accuracy or propagating errors due to using fewer bits to approximate floating point values. Empirically, we see that the quality of our results is not significantly impacted by using single-precision arithmetic. We also see consistent results across different GPU architectures: a Quadro RTX 8000, Titan RTX, and RTX A6000. Because PyHySCO is optimized to parallelize computations on GPU, the runtimes are faster on the GPUs compared to the Intel Xeon E5-4627 CPU.

\begin{table}[t]
\centering
\resizebox{\columnwidth}{!}{%
\begin{tabular}{||c|c||c|c||c|c||c|c||c|c||}
\hline
\multicolumn{2}{||c||}{} & \multicolumn{2}{c||}{RTX A6000} & \multicolumn{2}{c||}{Quadro RTX 8000} & \multicolumn{2}{c||}{Titan RTX} & \multicolumn{2}{c||}{Intel Xeon E5-4627}\\
\multicolumn{2}{||c||}{} & \multicolumn{2}{c||}{(GPU)} & \multicolumn{2}{c||}{(GPU)} & \multicolumn{2}{c||}{(GPU)} & \multicolumn{2}{c||}{(CPU)}\\
\cline{3-10}
\multicolumn{2}{||c||}{} & double & single & double & single & double & single & double & single \\
\hline\hline
\multirow{10}{*}{3T} & \multirow{2}{*}{Runtime (s)} & $12.7262$ & $9.5750$ & $13.4800$ & $7.8854$ & $13.1178$ & $7.5820$ & $34.3328$ & $27.1305$ \\
 & & $\pm 0.68$ & $\pm 0.58$ & $\pm 1.23$ & $\pm 0.91$ & $\pm 1.31$ & $\pm 0.98$ & $\pm 4.26$ & $\pm 3.09$ \\ \cline{2-10}
  & Optimization &$6.7947$ & $4.1562$ & $7.0133$ & $2.1065$ & $6.7862$ & $1.9327$ & $27.8682$ & $23.2062$ \\
 & Time (s) & $\pm 0.51$ & $\pm 0.39$ & $\pm 1.25$ & $\pm 0.90$ & $\pm 1.25$ & $\pm 0.62$ & $\pm 3.10$ & $\pm 3.07$ \\ \cline{2-10}
 & Relative & $82.7486$ & $82.7393$ & $82.7486$ & $82.7393$ & $82.7486$ & $82.7393$ & $82.486$ & $82.7393$ \\
 & Improvement & $\pm 3.49$ & $\pm 3.50$ & $\pm 3.49$ & $\pm 3.50$ & $\pm 3.49$ & $\pm 3.50$ & $\pm 3.49$ & $\pm 3.50$ \\ \cline{2-10}
 & \multirow{2}{*}{Loss Value} & $2.560e07$ & $2.562e07$ & $2.560e07$ & $2.562e07$ & $2.560e07$ & $2.562e07$ & $2.560e07$ & $2.562e07$ \\
 & & $\pm 7.66e06$ & $\pm 7.69e06$ & $\pm 7.66e06$ & $\pm 7.69e06$ & $\pm 7.66e06$ & $\pm 7.69e06$ & $\pm 7.66e06$ & $\pm 7.69e06$ \\ \cline{2-10}
 & Smoothness &  $3.848e04$ & $3.851e04$ & $3.848e04$ & $3.851e04$ & $3.848e04$ & $3.851e04$ & $3.848e04$ & $3.851e04$ \\
 & Reg. Value & $\pm 1.21e04$ & $\pm 1.21e04$ & $\pm 1.21e04$ & $\pm 1.21e04$ & $\pm 1.21e04$ & $\pm 1.21e04$ & $\pm 1.21e04$ & $\pm 1.21e04$ \\ \cline{2-10}
\hline\hline
\multirow{10}{*}{7T} & \multirow{2}{*}{Runtime (s)} & $17.2494$ & $11.9028$ & $21.0102$ & $9.3140$ & $18.6522$ & $9.4680$ & $82.1059$ & $33.4020$ \\
 & & $\pm 0.94$ & $\pm 0.44$ & $\pm 6.31$ & $\pm 0.99$ & $\pm 2.84$ & $\pm 2.19$ & $\pm 7.64$ & $\pm 4.15$ \\ \cline{2-10}
  & Optimization & $10.1298$ & $5.7460$ & $11.9937$ & $2.2579$ & $10.9617$ & $2.9775$ & $72.1380$ & $28.5530$ \\
 & Time (s) & $\pm 0.95$ & $\pm 0.42$ & $\pm 3.55$ & $\pm 0.64$ & $\pm 2.85$ & $\pm 2.11$ & $\pm 6.82$ & $\pm 4.09$ \\ \cline{2-10}
 & Relative & $85.7618$ & $85.7641$ & $85.7618$ & $85.7642$ & $85.7618$ & $85.7642$ & $85.7618$ & $85.7638$ \\
 & Improvement & $\pm 5.10$ & $\pm 5.10$ & $\pm 5.10$ & $\pm 5.10$ & $\pm 5.10$ & $\pm 5.10$ & $\pm 5.10$ & $\pm 5.10$ \\ \cline{2-10}
 & \multirow{2}{*}{Loss Value} & $4.143e07$ & $4.140e07$ & $4.143e07$ & $4.140e07$ & $4.143e07$ & $4.140e07$ & $4.1432e07$ & $4.1410e07$ \\
 & & $\pm 1.95e07$ & $\pm 1.95e07$ & $\pm 1.95e07$ & $\pm 1.95e07$ & $\pm 1.95e07$ & $\pm 1.95e07$ & $\pm 1.95e07$ & $\pm 1.95e07$ \\ \cline{2-10}
 & Smoothness & $5.634e04$ & $5.628e04$ & $5.634e04$ & $5.628e04$ & $5.634e04$ & $5.628e04$ & $5.634e04$ & $5.629e04$  \\
 & Reg. Value & $\pm 1.91e04$ & $\pm 1.91e04$ & $\pm 1.91e04$ & $\pm 1.91e04$ & $\pm 1.91e04$ & $\pm 1.91e04$ & $\pm 1.91e04$ & $\pm 1.91e04$ \\ \cline{2-10}
\hline\hline
\multirow{10}{*}{Sim.} & \multirow{2}{*}{Runtime (s)} & $106.92$ & $50.26$ & $127.38$ & $24.17$ & $125.25$ & $23.60$ & $851.18$ & $417.56$ \\
 & & $\pm 11.42$ & $\pm 3.52$ & $\pm 13.42$ & $\pm 1.31$ & $\pm 14.87$ & $\pm 3.78$ & $\pm 107.41$ & $\pm 55.33$ \\ \cline{2-10}
  & Optimization & $89.53$ & $35.53$ & $104.47$ & $7.93$ & $105.78$ & $9.32$ & $827.06$ & $402.04$ \\
 & Time (s) & $\pm 11.56$ & $\pm 4.13$ & $\pm 13.88$ & $\pm 0.92$ & $\pm 14.84$ & $\pm 3.81$ & $\pm 108.91$ & $\pm 56.58$ \\ \cline{2-10}
 & Relative & $76.27$ & $76.28$ & $76.27$ & $76.28$ & $76.27$ & $76.28$ & $76.27$ & $76.28$ \\
 & Improvement & $\pm 5.18$ & $\pm 5.18$ & $\pm 5.18$ & $\pm 5.18$ & $\pm 5.18$ & $\pm 5.18$ & $\pm 5.18$ & $\pm 5.18$ \\ \cline{2-10}
 & \multirow{2}{*}{Loss Value} & $6.07e07$ & $6.08e07$ & $6.07e07$ & $6.08e07$ & $6.07e07$ & $6.08e07$ & $6.07e07$ & $6.08e07$ \\
 & & $\pm 1.39e07$ & $\pm 1.40e07$ & $\pm 1.39e07$ & $\pm 1.40e07$ & $\pm 1.39e07$ & $\pm 1.40e07$ & $\pm 1.39e07$ & $\pm 1.40e07$ \\ \cline{2-10}
 & Smoothness & $9.53e04$ & $9.56e04$ & $9.53e04$ & $9.56e04$ & $9.53e04$ & $9.56e04$ & $9.54e04$ & $9.56e04$ \\
 & Reg. Value & $\pm 2.71e04$ & $\pm 2.72e04$ & $\pm 2.71e04$ & $\pm 2.73e04$ & $\pm 2.71e04$ & $\pm 2.73e04$ & $\pm 2.71e04$ & $\pm 2.73e04$ \\ \cline{2-10}
\hline\hline
\end{tabular}
}%
\caption{The speed and quality of PyHySCO optimization with Gauss Newton on three different GPUs and a CPU in both single (float 32) and double (float 64) precision arithmetic. The relative improvement, loss value, and smoothness value are evaluated in double precision in all cases. Results are shown for both 3T and 7T data from the Human Connectome Project \cite{HCPdata} and simulated data. There is a great speedup when calculating in single precision without losing the quality of correction, and the speedup of PyHySCO using a GPU is clear compared to the CPU.}
\label{table:precision}
\end{table}

\subsection{Comparison of PyHySCO with HySCO and TOPUP}
\label{subs:HyscoTopup}
\begin{figure}
    \centering
    \input{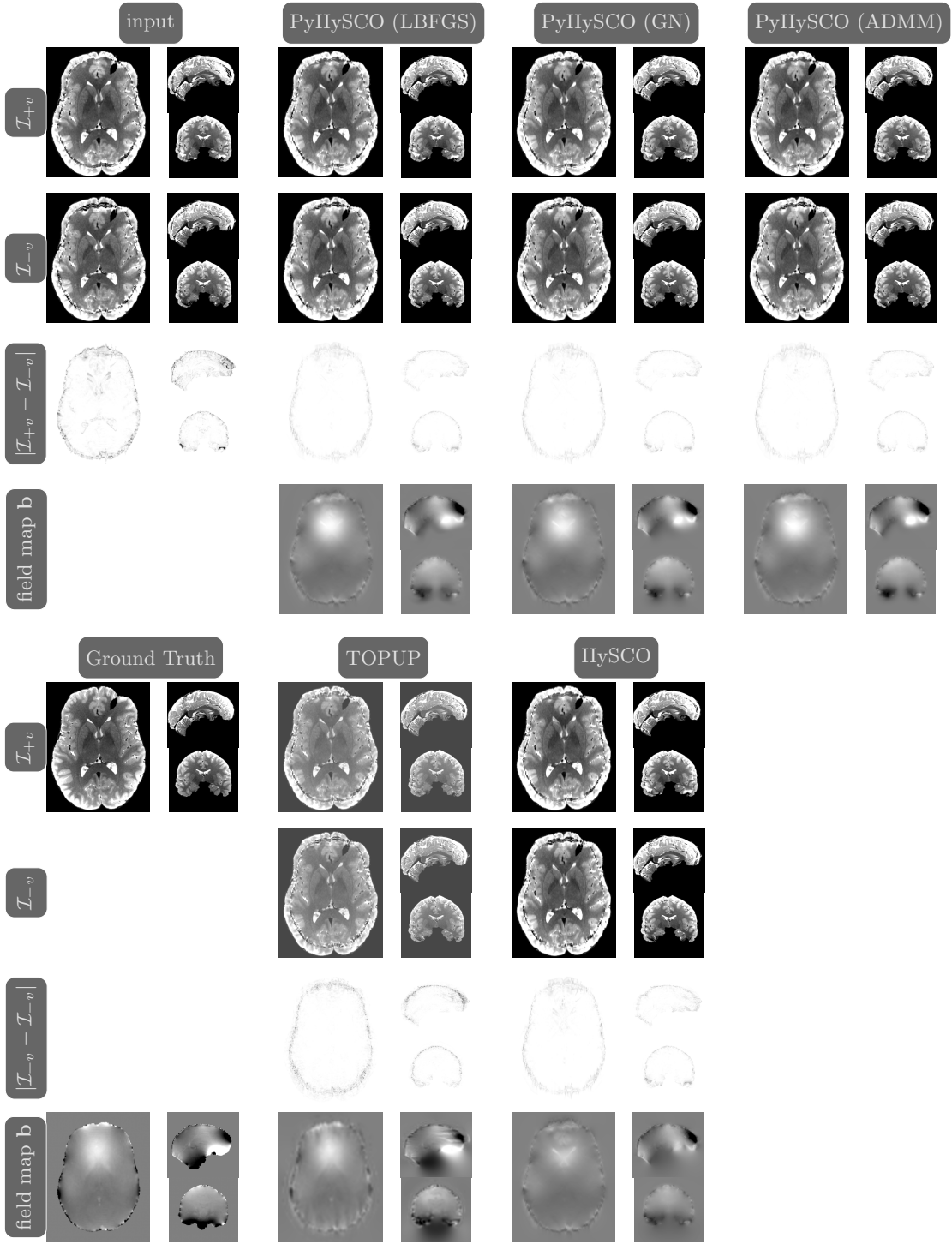}
    \caption{Visualization of resulting field maps and images for one subject from the simulated dataset (Subject ID 105014). The first column in the top half shows the input data and the first column on the bottom shows the ground truth T2w image and field map. The remaining columns show the results from PyHySCO using LBFGS, GN-PCG, and ADMM, TOPUP, and HySCO. For each optimization, the top two rows are the pair of images with opposite phase encoding directions, and the third row shows the absolute difference (with inverted color) between the pair of images. The bottom row shows the field maps estimated for each method.  PyHySCO achieves similar image distance and field map smoothness improvements in less computational time.}
    \label{fig:resultssim}
\end{figure}

Table~\ref{table:optimizers} reports the runtime and correction quality for PyHySCO using GN-PCG, HySCO, and TOPUP. On real 3T and 7T data, PyHySCO achieves lower loss and higher relative improvement between corrected images than HySCO, and higher relative improvement than TOPUP. The runtime on CPU for real data is 1-2 minutes for HySCO and over 1 hour for TOPUP, while PyHySCO on GPU has runtimes of 10-13 seconds. For the simulated dataset, PyHySCO requires an average of 1 minute on GPU, HySCO an average of 12.6 minutes on CPU, and TOPUP an average of 8.5 hours on CPU. Using the ground truth field maps from the simulated dataset, PyHySCO achieves the lowest average field map relative error, 14.48\%, compared to 19.70\% for HySCO and 16.36\% for TOPUP. Figures \ref{fig:results3T}, \ref{fig:results7T}, and \ref{fig:resultssim} show the field map and corrected images for one example subject from each dataset. The results of the methods are similar, and the resulting field maps are comparable to those of the existing tools, HySCO and TOPUP, while PyHySCO is considerably faster.

\begin{table}[t]
\begin{center}
\begin{tabular}{||c|c||c||c||c||}
\hline
\multicolumn{2}{||c||}{} & \multicolumn{1}{c||}{PyHySCO} & \multicolumn{1}{c||}{HySCO} & \multicolumn{1}{c||}{TOPUP}\\
\cline{3-5}
\hline\hline
\multirow{8}{*}{3T} & \multirow{2}{*}{Runtime (s)} & $10.37$ & $65.06$ & $4022.56$ \\
 & & $\pm 0.87$ & $\pm 8.64$ & $\pm 73.11$ \\ \cline{2-5}
 & Relative & $82.74$ & $78.98$ & $54.36$ \\
 & Improvement & $\pm 3.50$ & $\pm 6.39$ & $\pm 17.08$ \\ \cline{2-5}
 & \multirow{2}{*}{Loss Value} & $2.56e07$ & $4.13e07$ & \multirow{2}{*}{N/A} \\
 & & $\pm 7.69e06$ & $\pm 1.38e07$ & \\ \cline{2-5}
 & Smoothness & $3.85e04$ & $7.84e04$ & \multirow{2}{*}{N/A} \\
 & Reg. Value &  $\pm 1.21e04$ & $\pm 3.01e04$ &  \\ \cline{2-5}
\hline\hline
\multirow{8}{*}{7T} & \multirow{2}{*}{Runtime (s)} & $13.62$ & $120.92$ & $3713.51$ \\
 & & $\pm 2.38$ & $\pm 19.61$ & $\pm 63.04$ \\ \cline{2-5}
 & Relative & $85.76$ & $80.43$ & $74.51$ \\
 & Improvement &  $\pm 5.10$ & $\pm 10.46$ & $\pm 9.13$ \\ \cline{2-5}
 & \multirow{2}{*}{Loss Value} & $4.14e07$ & $5.87e07$ & \multirow{2}{*}{N/A} \\
 & & $\pm 1.95e07$ & $\pm 2.48e07$ &  \\ \cline{2-5}
 & Smoothness &  $5.63e04$ & $8.03e04$ & \multirow{2}{*}{N/A} \\
 & Reg. Value & $\pm 1.91e04$ & $\pm 3.68e04$ &  \\ \cline{2-5}
\hline\hline
\multirow{8}{*}{Simulated} & \multirow{2}{*}{Runtime (s)} &$55.26$ & $757.65$ & $30854.18$ \\
 & & $\pm 3.86$ & $\pm 96.26$ & $\pm 568.11$ \\ \cline{2-5}
 & Relative &  $76.28$ & $69.53$ & $17.56$ \\
 & Improvement & $\pm 5.18$ & $\pm 5.10$ & $\pm 28.14$ \\ \cline{2-5}
 & Loss & $6.08e07$ & $6.07e07$ & \multirow{2}{*}{N/A} \\
 & Value &  $\pm1.40e07$ & $\pm 1.51e07$ &  \\ \cline{2-5}
 & Smoothness & $9.56e04$ & $6.10e04$ & \multirow{2}{*}{N/A} \\
 & Reg. Value & $\pm 2.72e04$ & $\pm 1.60e04$ &  \\ \cline{2-5}
  & Relative Error & $14.48$ & $19.70$ & $16.37$ \\
  & (Field Map) & $\pm 7.71$ & $\pm 11.70$ & $\pm 3.60$ \\ \cline{2-5}
\hline\hline
\end{tabular}
\caption{The speed and quality of optimization for TOPUP, HySCO, and PyHySCO. PyHySCO uses Gauss Newton and optimizes in single precision on GPU. HySCO and TOPUP optimize on CPU using the default configurations. Results are reported for 3T and 7T data from the Human Connectome Project \cite{HCPdata} and the simulated distortion data.}
\label{table:optimizers}
\end{center}
\end{table}

\begin{figure}
    \centering
    \input{Fig5aCompareResults3T}
    \caption{Visualization of resulting field maps and images for one subject from the 3T dataset (Subject ID 211619). The first column in the top half shows the input data. The remaining columns show the results from PyHySCO using LBFGS, GN-PCG, and ADMM compared with TOPUP and HySCO. For each optimization, the top two rows are the pair of images with opposite phase encoding directions, and the third row shows the absolute difference (with inverted color) between the pair of images. The bottom row shows the field maps estimated for each method. PyHySCO achieves similar image distance and field map smoothness improvements in less computational time. }
    \label{fig:results3T}
\end{figure}

\begin{figure}
    \centering
    \input{Fig5bCompareResults7T}
    \caption{Visualization of resulting field maps and images for one subject from the 7T dataset (Subject ID 825048). The first column in the top half shows the input data. The remaining columns show the results from PyHySCO using LBFGS, GN-PCG, and ADMM compared with TOPUP and HySCO. For each optimization, the top two rows are the pair of images with opposite phase encoding directions, and the third row shows the absolute difference (with inverted color) between the pair of images. The bottom row shows the field maps estimated for each method. PyHySCO achieves similar image distance and field map smoothness improvements in less computational time.}
    \label{fig:results7T}
\end{figure}

\section{Discussion}\label{sec:4discussion}
The PyHySCO toolbox accurately and robustly corrects susceptibility artifacts in EPIs acquired using the Reverse Gradient Polarity acquisition. 
In numerous experiments with real and simulated data, it achieves similar correction quality to the leading RGP toolboxes TOPUP and HySCO while having a time-to-solution in the order of timings reported for pre-trained deep learning approaches.
Compared to the latter class of methods, it is important to highlight that PyHySCO does not require any training and is based on a physical distortion model, which helps generalize to different scanners, image acquisition parameters, and anatomies. 

PyHySCO's modular design invites improvements and contributions. 
The toolbox is based on PyTorch, which provides hardware support and other functionality, including automatic differentiation.
In our experiments,  correction quality is hardware and precision-independent, but a considerable speedup is realized on GPUs with single precision (32-bit) arithmetic. 
The reduced computational time is mostly attributed to the effective use of multithreading and parallelism on modern hardware. 

PyHySCO introduces a novel, fast, and effective optimal transport-based initialization scheme. The initial estimate of the field map already substantially reduces the distance between the images with opposite phase encoding directions. In our experiments, the non-smoothness of the initial field map can be corrected by applying a Gaussian blur and a few optimization steps to the full image resolution.  

PyHySCO's three optimization algorithms achieve comparable correction results but have different computational costs. The ADMM algorithm takes advantage of the separable structure of the optimization problem to enhance parallelism but requires more iterations than GN-PCG. 
While this results in longer runtimes in our examples, the method could be more scalable for datasets of considerably higher resolution. 
For the relatively standard image sizes of about $200 \times 200 \times 132$, the default GN-PCG algorithm is most effective.
Both customized optimization algorithms are more efficient than our comparison, LBFGS.

PyHySCO can be interfaced directly in Python or run in a batch mode via the command line. The latter makes it a drop-in replacement for other RGP tools in MRI post-processing pipelines.

\section{Conclusion}\label{sec:conclusion}
PyHySCO offers RGP-based correction with high accuracy at the cost similar to pre-trained learning-based methods. Our implementation is based on PyTorch and makes efficient use of modern hardware accelerators such as GPUs. We show the accuracy and efficiency of PyHySCO on real and simulated three-dimensional volumes of various field strengths and phase encoding axes. Our results show that PyHySCO achieves a correction of comparable quality to leading physics-based methods in a fraction of the time.

\subsection*{Data and Code Availability}
\noindent The source code, examples, and documentation for PyHySCO are available at the following repository: \hyperlink{https://github.com/EmoryMLIP/PyHySCO}{https://github.com/EmoryMLIP/PyHySCO}. The Python package for PyHySCO can be installed via pip and be downloaded from: \hyperlink{https://pypi.org/project/PyHySCO/}{https://pypi.org/project/PyHySCO/}.

Data were provided by the Human Connectome Project, WU-Minn Consortium (Principal Investigators: David Van Essen and Kamil Ugurbil; 1U54MH091657) funded by the 16 NIH Institutes and Centers that support the NIH Blueprint for Neuroscience Research; and by the McDonnell Center for Systems Neuroscience at Washington University.

\subsection*{Author Contributions}
AJ and LR contributed to the code development and preparation of this manuscript. AJ ran the experiments and examples for the manuscript.

\subsection*{Acknowledgements}
AJ is supported by the National Science Foundation Graduate Research Fellowship under Grant No. 1937971. The work was also supported in part by the NSF awards DMS 1751636 and DMS 2038118.

\subsection*{Declaration of Competing Interests}
The authors declare no competing interests.

\bibliographystyle{apalike}
\bibliography{references}

\begin{thebibliography}{}

\bibitem[Alkilani et~al., 2023]{alkilani2023fd}
Alkilani, A.~Z., {\c{C}}ukur, T., and Saritas, E.~U. (2023).
\newblock Fd-net: An unsupervised deep forward-distortion model for susceptibility artifact correction in epi.
\newblock {\em arXiv preprint arXiv:2303.10436}.

\bibitem[Andersson et~al., 2003]{AnderssonEtAl2003}
Andersson, J. L.~R., Skare, S., and Ashburner, J. (2003).
\newblock {How to correct susceptibility distortions in spin-echo echo-planar images: application to diffusion tensor imaging}.
\newblock {\em NeuroImage}, 20(2):870--888.

\bibitem[Antun et~al., 2020]{AntunEtAl2020}
Antun, V., Renna, F., Poon, C., Adcock, B., and Hansen, A.~C. (2020).
\newblock {On instabilities of deep learning in image reconstruction and the potential costs of AI}.
\newblock {\em Proceedings of the National Academy of Sciences}, 117(48):30088--30095.

\bibitem[Bowtell et~al., 1994]{bowtell1994correction}
Bowtell, R., McIntyre, D., Commandre, M., Glover, P., and Mansfield, P. (1994).
\newblock Correction of geometric distortion in echo planar images.
\newblock In {\em Soc. Magn. Res. Abstr}, volume~2, page 411.

\bibitem[Boyd et~al., 2011]{boyd2011admm}
Boyd, S., Parikh, N., Chu, E., Peleato, B., Eckstein, J., et~al. (2011).
\newblock Distributed optimization and statistical learning via the alternating direction method of multipliers.
\newblock {\em Foundations and Trends{\textregistered} in Machine learning}, 3(1):1--122.

\bibitem[Cai et~al., 2021]{cai2021prequal}
Cai, L.~Y., Yang, Q., Hansen, C.~B., Nath, V., Ramadass, K., Johnson, G.~W., Conrad, B.~N., Boyd, B.~D., Begnoche, J.~P., Beason-Held, L.~L., et~al. (2021).
\newblock Prequal: An automated pipeline for integrated preprocessing and quality assurance of diffusion weighted mri images.
\newblock {\em Magnetic resonance in medicine}, 86(1):456--470.

\bibitem[Chang and Fitzpatrick, 1992]{ChangFitzpatrick1992}
Chang, H. and Fitzpatrick, J.~M. (1992).
\newblock {A Technique for Accurate Magnetic-Resonance-Imaging in the Presence of Field Inhomogeneities}.
\newblock {\em Medical Imaging, IEEE Transactions on}, 11(3):319--329.

\bibitem[Chen et~al., 2022]{chen2022DLReview}
Chen, Z., Pawar, K., Ekanayake, M., Pain, C., Zhong, S., and Egan, G.~F. (2022).
\newblock Deep learning for image enhancement and correction in magnetic resonance imaging—state-of-the-art and challenges.
\newblock {\em Journal of Digital Imaging}, pages 1--27.

\bibitem[Cooley et~al., 1969]{cooley1969convfft}
Cooley, J.~W., Lewis, P.~A., and Welch, P.~D. (1969).
\newblock The fast fourier transform and its applications.
\newblock {\em IEEE Transactions on Education}, 12(1):27--34.

\bibitem[D{\'a}vid et~al., 2024]{Dvid2024ACIDAC}
D{\'a}vid, G., Fricke, B., Oeschger, J.~M., Ruthotto, L., Fritz, F.~J., Ohana, O., Sauvigny, T., Freund, P., Tabelow, K., and Mohammadi, S. (2024).
\newblock Acid: A comprehensive toolbox for image processing and modeling of brain, spinal cord, and ex vivo diffusion mri data.
\newblock {\em bioRxiv}.

\bibitem[Duong et~al., 2020a]{duong2020susceptibility}
Duong, S., Phung, S.~L., Bouzerdoum, A., Taylor, H.~B., Puckett, A., and Schira, M.~M. (2020a).
\newblock Susceptibility artifact correction for sub-millimeter fmri using inverse phase encoding registration and t1 weighted regularization.
\newblock {\em Journal of Neuroscience Methods}, 336:108625.

\bibitem[Duong et~al., 2020b]{duong2020unsupervised}
Duong, S.~T., Phung, S.~L., Bouzerdoum, A., and Schira, M.~M. (2020b).
\newblock An unsupervised deep learning technique for susceptibility artifact correction in reversed phase-encoding epi images.
\newblock {\em Magnetic Resonance Imaging}, 71:1--10.

\bibitem[Duong et~al., 2021]{Duong2021}
Duong, S. T.~M., Phung, S.~L., Bouzerdoum, A., Ang, S.~P., and Schira, M.~M. (2021).
\newblock Correcting susceptibility artifacts of mri sensors in brain scanning: A 3d anatomy-guided deep learning approach.
\newblock {\em Sensors}, 21(7).

\bibitem[Esteban et~al., 2014]{esteban2014simulation}
Esteban, O., Daducci, A., Caruyer, E., O'Brien, K., Ledesma-Carbayo, M.~J., Bach-Cuadra, M., and Santos, A. (2014).
\newblock Simulation-based evaluation of susceptibility distortion correction methods in diffusion mri for connectivity analysis.
\newblock In {\em 2014 IEEE 11th International Symposium on Biomedical Imaging (ISBI)}, pages 738--741. IEEE.

\bibitem[Graham et~al., 2017]{graham2017quantitative}
Graham, M.~S., Drobnjak, I., Jenkinson, M., and Zhang, H. (2017).
\newblock Quantitative assessment of the susceptibility artefact and its interaction with motion in diffusion mri.
\newblock {\em PloS one}, 12(10):e0185647.

\bibitem[Gu and Eklund, 2019]{gu2019evaluation}
Gu, X. and Eklund, A. (2019).
\newblock Evaluation of six phase encoding based susceptibility distortion correction methods for diffusion mri.
\newblock {\em Frontiers in neuroinformatics}, 13:76.

\bibitem[Hansen et~al., 2006]{hansen2006deblurring}
Hansen, P., Nagy, J., and O'Leary, D. (2006).
\newblock {\em Deblurring Images: Matrices, Spectra, and Filtering}.
\newblock Fundamentals of Algorithms. Society for Industrial and Applied Mathematics.

\bibitem[Hestenes and Stiefel, 1952]{HestenesStiefel1952}
Hestenes, M.~R. and Stiefel, E. (1952).
\newblock {Methods of Conjugate Gradients for Solving Linear Systems}.
\newblock {\em Journal of Research of the National Bureau of Standards}, 49(6):409 436.

\bibitem[Holland et~al., 2010]{HollandEtAl2010}
Holland, D., Kuperman, J.~M., and Dale, A.~M. (2010).
\newblock {Efficient correction of inhomogeneous static magnetic field-induced distortion in Echo Planar Imaging}.
\newblock {\em NeuroImage}, 50(1):175--183.

\bibitem[Hu et~al., 2020]{Hu2020}
Hu, Z., Wang, Y., Zhang, Z., Zhang, J., Zhang, H., Guo, C., Sun, Y., and Guo, H. (2020).
\newblock Distortion correction of single-shot epi enabled by deep-learning.
\newblock {\em NeuroImage}, 221:117--170.

\bibitem[Irfanoglu et~al., 2015]{irfanoglu2015dr}
Irfanoglu, M.~O., Modi, P., Nayak, A., Hutchinson, E.~B., Sarlls, J., and Pierpaoli, C. (2015).
\newblock Dr-buddi (diffeomorphic registration for blip-up blip-down diffusion imaging) method for correcting echo planar imaging distortions.
\newblock {\em Neuroimage}, 106:284--299.

\bibitem[Liu and Nocedal, 1989]{lbfgs}
Liu, D.~C. and Nocedal, J. (1989).
\newblock On the limited memory bfgs method for large scale optimization.
\newblock {\em Mathematical programming}, 45(1):503--528.

\bibitem[Macdonald and Ruthotto, 2017]{MacdonaldRuthotto2017}
Macdonald, J. and Ruthotto, L. (2017).
\newblock {Improved Susceptibility Artifact Correction of Echo Planar MRI using the Alternating Direction Method of Multipliers}.
\newblock {\em Journal of Mathematical Imaging and Vision}, 60(2):268--282.

\bibitem[Modersitzki, 2009]{Modersitzki2009}
Modersitzki, J. (2009).
\newblock {\em {FAIR: flexible algorithms for image registration}}, volume~6 of {\em Society for Industrial and Applied Mathematics (SIAM), Philadelphia, PA}.
\newblock Society for Industrial and Applied Mathematics (SIAM), Philadelphia, PA.

\bibitem[Nocedal and Wright, 1999]{nocedal1999numerical}
Nocedal, J. and Wright, S.~J. (1999).
\newblock {\em Numerical optimization}.
\newblock Springer.

\bibitem[Paszke et~al., 2019]{pytorch}
Paszke, A., Gross, S., Massa, F., Lerer, A., Bradbury, J., Chanan, G., Killeen, T., Lin, Z., Gimelshein, N., Antiga, L., Desmaison, A., Kopf, A., Yang, E., DeVito, Z., Raison, M., Tejani, A., Chilamkurthy, S., Steiner, B., Fang, L., Bai, J., and Chintala, S. (2019).
\newblock Pytorch: An imperative style, high-performance deep learning library.
\newblock In Wallach, H., Larochelle, H., Beygelzimer, A., d\textquotesingle Alch\'{e}-Buc, F., Fox, E., and Garnett, R., editors, {\em Advances in Neural Information Processing Systems}, volume~32. Curran Associates, Inc.

\bibitem[Penny et~al., 2007]{SPM2007}
Penny, W.~D., Friston, K.~J., Ashburner, J.~T., Kiebel, S.~J., and Nichols, T.~E. (2007).
\newblock {\em Statistical parametric mapping: the analysis of functional brain images}.
\newblock Elsevier.

\bibitem[Peyr{\'e} et~al., 2017]{peyreOT}
Peyr{\'e}, G., Cuturi, M., et~al. (2017).
\newblock Computational optimal transport.
\newblock {\em Center for Research in Economics and Statistics Working Papers}, 2017-86.

\bibitem[Ruthotto et~al., 2012]{RuthottoEtAl2012EPI}
Ruthotto, L., Kugel, H., Olesch, J., Fischer, B., Modersitzki, J., Burger, M., and Wolters, C.~H. (2012).
\newblock {Diffeomorphic susceptibility artifact correction of diffusion-weighted magnetic resonance images}.
\newblock {\em Physics in Medicine and Biology}, 57(18):5715--5731.

\bibitem[Ruthotto et~al., 2013]{RuthottoEtAl2013hysco}
Ruthotto, L., Mohammadi, S., Heck, C., Modersitzki, J., and Weiskopf, N. (2013).
\newblock {Hyperelastic Susceptibility Artifact Correction of DTI in SPM}.
\newblock In {\em Bildverarbeitung fuer die Medizin}, pages 344--349, Berlin, Heidelberg. Springer, Berlin, Heidelberg.

\bibitem[Saad, 2003]{Saad2003}
Saad, Y. (2003).
\newblock {\em {Iterative Methods for Sparse Linear Systems}}.
\newblock SIAM. SIAM.

\bibitem[Smith et~al., 2004]{smith2004FSL}
Smith, S.~M., Jenkinson, M., Woolrich, M.~W., Beckmann, C.~F., Behrens, T.~E., Johansen-Berg, H., Bannister, P.~R., De~Luca, M., Drobnjak, I., Flitney, D.~E., et~al. (2004).
\newblock Advances in functional and structural mr image analysis and implementation as fsl.
\newblock {\em Neuroimage}, 23:S208--S219.

\bibitem[Snoussi et~al., 2021]{Snoussi2021_spinalcordeval}
Snoussi, H., Cohen-Adad, J., Commowick, O., Combes, B., Bannier, E., Leguy, S., Kerbrat, A., Barillot, C., and Caruyer, E. (2021).
\newblock Evaluation of distortion correction methods in diffusion mri of the spinal cord.

\bibitem[Stehling et~al., 1991]{Stehling1991EPI}
Stehling, M.~K., Turner, R., and Mansfield, P. (1991).
\newblock Echo-planar imaging: Magnetic resonance imaging in a fraction of a second.
\newblock {\em Science}, 254(5028):43--50.

\bibitem[Tax et~al., 2022]{tax2022s}
Tax, C.~M., Bastiani, M., Veraart, J., Garyfallidis, E., and Irfanoglu, M.~O. (2022).
\newblock What’s new and what’s next in diffusion mri preprocessing.
\newblock {\em NeuroImage}, 249:118830.

\bibitem[Van~Essen et~al., 2012]{HCPdata}
Van~Essen, D.~C., Ugurbil, K., Auerbach, E., Barch, D., Behrens, T.~E., Bucholz, R., Chang, A., Chen, L., Corbetta, M., Curtiss, S.~W., et~al. (2012).
\newblock The human connectome project: a data acquisition perspective.
\newblock {\em Neuroimage}, 62(4):2222--2231.

\bibitem[Wu et~al., 2008]{wu2008comparison}
Wu, M., Chang, L.-C., Walker, L., Lemaitre, H., Barnett, A.~S., Marenco, S., and Pierpaoli, C. (2008).
\newblock Comparison of epi distortion correction methods in diffusion tensor mri using a novel framework.
\newblock In {\em International Conference on Medical Image Computing and Computer-Assisted Intervention}, pages 321--329. Springer.

\bibitem[Zahneisen et~al., 2020]{zahneisen2020deep}
Zahneisen, B., Baeumler, K., Zaharchuk, G., Fleischmann, D., and Zeineh, M. (2020).
\newblock Deep flow-net for epi distortion estimation.
\newblock {\em Neuroimage}, 217:116886.

\end{thebibliography}

\end{document}